\documentclass[10pt]{amsart}

\usepackage{latexsym,amssymb,mathrsfs,amsmath,amsthm,amscd,psfrag,epic,eepic,epsfig}
\usepackage[dvips]{color}

\newtheorem{theorem}{Theorem}
\newtheorem{lemma}[theorem]{Lemma}
\newtheorem{prop}[theorem]{Proposition}
\newtheorem{cor}[theorem]{Corollary}

\theoremstyle{definition}
\newtheorem{definition}[theorem]{Definition}

\theoremstyle{definition}
\newtheorem*{remark}{Remark}
\newtheorem*{remarks}{Remarks}
\newtheorem*{convention}{Convention}

\theoremstyle{remark}
\newtheorem*{acknowledgement}{Acknowledgement}

\theoremstyle{plain}

\theoremstyle{plain}

\newenvironment{claim}[1]{\smallskip\\{\it Claim #1.}}{}
\newenvironment{pc}{\smallskip\\}{}
\newenvironment{step}[1]{\smallskip\\{\it Step #1.}}{}

\newcommand{\Z}{\mathbb{Z}}
\newcommand{\R}{\mathbb{R}}
\newcommand{\N}{\mathbb{N}}
\newcommand{\Q}{\mathbb{Q}}
\newcommand{\C}{\mathbb{C}}

\newcommand{\Hyp}{\mathbb{H}}

\newcommand{\e}{\varepsilon}

\newcommand{\action}{\mathcal{A}}
\newcommand{\M}{\mathcal{M}}
\newcommand{\de}{\mathrm{d}}
\newcommand{\De}{\mathrm{D}}

\newcommand{\Ce}{\tilde{C}}
\newcommand{\Nn}{\tilde{N}}
\newcommand{\HF}{HF_*}
\newcommand{\CF}{CF_*}

\newcommand{\CM}{CM_*}

\newcommand{\CoM}{CM^*}

\newcommand{\Tc}{\mathcal{T}}
\newcommand{\tc}{\mathfrak{t}}
\newcommand{\Ac}{\mathcal{A}}
\newcommand{\Bc}{\mathcal{B}}
\newcommand{\Ec}{\mathcal{E}}
\newcommand{\Vc}{\mathcal{V}}
\newcommand{\Cc}{\mathcal{C}}

\newcommand{\I}{[0,1]}

\renewcommand{\to}{\rightarrow}
\newcommand{\To}{\longrightarrow}
\newcommand{\Mapsto}{\longmapsto}
\newcommand{\Iff}{\Longleftrightarrow}

\newcommand{\then}{\Rightarrow}

\newcommand{\inclusion}{\hookrightarrow}

\newcommand{\p}{\partial}

\DeclareMathOperator{\id}{id}

\DeclareMathOperator{\im}{im}

\DeclareMathOperator{\mind}{ind}
\DeclareMathOperator{\pr}{pr}

\DeclareMathOperator{\area}{area}
\DeclareMathOperator{\crit}{Crit}
\DeclareMathOperator{\fix}{Fix}

\DeclareMathOperator{\symp}{Symp}
\DeclareMathOperator{\diff}{Diff}

\DeclareMathOperator{\sign}{sign}
\DeclareMathOperator{\flux}{Flux}
\DeclareMathOperator{\supp}{supp}

\newcommand{\tf}{T_{\phi}}

\newcommand{\Of}{\Omega_{\phi}}
\newcommand{\ophi}{\omega_{\phi}}

\title[Floer homology of algebraically finite mapping classes]
{Floer homology of algebraically finite \\ mapping classes}
\date{\today}
\author{Ralf Gautschi}
\thanks
{Research partially supported by the SNF (project number 21-64937.01)}
\address{ETH Z\"urich, HG G36.1, R\"amistrasse 101, 8092 Z\"urich}
\email{ralf.gautschi@math.ethz.ch}

\begin{document}

\begin{abstract}
Using symplectic Floer homology, Seidel associated a module to each mapping
class of a compact connected oriented two-manifold of genus bigger than one.
We compute this module for mapping classes which do not have any pseudo-Anosov
components in the sense of Thurston's theory of surface diffeomorphisms. 
The Nielsen-Thurston representative of such a class is shown
to be monotone. The formula for the Floer homology is obtained from 
a topological separation of fixed points and a separation
mechanism for Floer connecting orbits. As examples, we consider the geometric
monodromy of isolated plane curve singularities. In this case, the formula for
the Floer homology is particularly simple.  
\end{abstract}

\maketitle

\section{Introduction}

This work is concerned with the computation of Floer homology groups of
symplectomorphisms of compact 2-manifolds.  
In the case of a 2-sphere, each symplectomorphism is exact and its Floer
homology equals, by Floer's original work \cite{F4}, the singular homology of
the 2-sphere.  
In the case of a torus, the Floer homology of linear symplectomorphisms was
computed by Po{\'z}niak~\cite{Po}. 

In this article, we consider the case of a compact connected oriented 2-manifold
$\Sigma$ of genus $\geq2$.
It was shown by Seidel~\cite{S2}, that to every mapping class $g$ of $\Sigma$,
there is associated a $\Z_2$-graded vector space $\HF(g)$ over the field
$\Z_2$, which is equipped with
an additional module structure over $H^*(\Sigma;\Z_2)$. 
To put it short, $\HF(g)$ is the symplectic Floer
homology of an area preserving and so-called monotone representative of $g$.

The first computational result in this context is also due to
Seidel. 
In \cite{S1} it was shown, that
if $g$ is the mapping class of a positive Dehn twist along an embedded circle
$C\subset\Sigma$, then $\HF(g)\cong H_*(\Sigma\setminus C;\Z_2)$. 

Our starting point is the following definition given in detail in 
Section~\ref{se:diff}: an orientation preserving diffeomorphism $\phi$ of
$\Sigma$ is called of finite type, if $\Sigma$ can be obtained by piecing together
$\phi$-invariant 2-manifolds $\Sigma'$ with boundary such that $\phi|\Sigma'$
is either periodic, a flip-twist map or a twist map without fixed
points. For the terminology we refer to Section~\ref{se:diff}. 
The significance of this definition lies in the following fact: due to 
Nielsen-Thurston theory of surface diffeomorphisms, every mapping class
without pseudo-Anosov components admits a representative which is of finite
type. Such a mapping class is called algebraically finite.

Our main result is a formula for the Floer homology of finite type diffeomorphisms.
We use the following notation. 
If $\phi$ is such a diffeomorphism, we denote by $\Sigma_0$ the
union of the $\Sigma'$ where $\phi$ restricts to the identity.
By $\p_+\Sigma_0$, we denote the union of components of $\p\Sigma_0$
where in a neighborhood, $\phi$ is a right-handed twist.
\begin{theorem}\label{thm:main1}
If $\phi$ is a diffeomorphism of finite type, then $\phi$ is monotone with
respect to some $\phi$-invariant area form and 
$$
\HF(\phi) \cong
H_*(\Sigma_0,\p_+\Sigma_0;\Z_2)\oplus
\Z_2^{\Lambda(\phi|\Sigma\setminus\Sigma_0)}.
$$
Here, $\Lambda$ denotes the Lefschetz number. 
The $H^*(\Sigma;\Z_2)$-action on the r.h.s. is split and given as follows.
On the first summand, it is the ordinary cap product.
On the second, $1\in H^0(\Sigma;\Z_2)$ acts by 
the identity and any element of $H^1(\Sigma;\Z_2)\oplus H^2(\Sigma;\Z_2)$ by
the zero map.  
\end{theorem}

As suggested by the formula above, the Floer complex of $\phi$ splits into a
complex associated to $\phi|\Sigma_0$ and one associated to
$\phi|\Sigma\setminus\Sigma_0$. On one hand this follows from a purely
topological separation of fixed points for finite type diffeomorphisms. 
On the other hand, there is also a separation mechanism for Floer connecting
orbits, i.e. after suitably perturbing $\phi$, every connecting orbit
starting and ending in $\Sigma_0$ does not leave $\Sigma_0$. For the precise
statement we refer to Section~\ref{se:corbits}.

A natural source of examples is provided by the theory of singularities. 
To every isolated plane curve singularity is associated a compact
connected oriented 2-manifold with boundary, the Milnor fiber, and an isotopy
class of 
orientation preserving diffeomorphisms of the Milnor fiber which are
the identity near the boundary, called the geometric monodromy.  
The precise definitions are given in Section~\ref{se:singularity}.
\begin{theorem}\label{thm:main2}
Let $M\subset\Sigma$ be the Milnor fiber of an isolated plane curve
singularity and $g$ be the mapping class which is obtained by extending the
geometric monodromy trivially to $\Sigma$. Then 
$$
\HF(g) \cong H_*(\Sigma,M;\Z_2),
$$
where $H^*(\Sigma;\Z_2)$ acts by cap product.
\end{theorem}
A special case of Theorem~\ref{thm:main2} is the following
generalization of Seidel's formula for the Floer homology of a Dehn
twist \cite{S1}. For the notation see Section~\ref{se:singularity}.
\begin{cor}\label{cor:dehn1}
Let $(C_1,\dots,C_k)$ be an $A_k$-configuration of circles in $\Sigma$. 
Let $g$ be the mapping class of the product 
$\tau_{1}\circ\cdots\circ\tau_{k}$, 
where $\tau_i$ denotes the positive Dehn twist along $C_i$. Then 
$$
\HF(g) \cong H_*(\Sigma,C_1\cup\cdots\cup C_k;\Z_2),
$$
where $H^*(\Sigma;\Z_2)$ acts by cap product. The same formula holds, 
if $g$ is the class of $\tau_{\sigma 1}\circ\cdots\circ\tau_{\sigma k}$, where
$\sigma$ is a cyclic permutation of $k$ elements.
\end{cor}
The proof of Theorem~\ref{thm:main2} is given in
Section~\ref{se:singularity}. It relies, besides on Theorem~\ref{thm:main1},
on an additional result about the geometric monodromy.
This result is stated in Proposition~\ref{prop:monodromy} and proven in
Appendix~\ref{ap:monodromy}. 
The proof uses the well known theorem of A'Campo~\cite{AC2} and
L{\^e}~\cite{L} that the Lefschetz number of the geometric monodromy
vanishes, together with the theory of splice diagrams which is due to Eisenbud
and Neumann~\cite{EN}. 
Appendix~B includes a summary of the relevant results on splice diagrams from
\cite{EN} and is therefore quite lengthy.
We would like to point out, that Proposition~\ref{prop:monodromy}
also follows from A'Campo's work~\cite{AC1}, \cite{AC4}.

Finally, we mention that there is a version of
Floer homology theory for diffeomorphisms of compact oriented 2-manifolds
with boundary.  
It assigns a pair of vector spaces $\HF(g,+),\HF(g,-)$ to every isotopy
class $g$ of orientation preserving diffeomorphisms which are the identity
near the boundary. The sign corresponds to two different 
perturbations of $g$ near the boundary. In Appendix \ref{ap:open} we
give an outline on this version of Floer homology.  
Using Proposition~\ref{prop:monodromy}, we confirm in
Section~\ref{se:singularity} a conjecture of
Seidel \cite[Page 23]{S5} in the case of plane curve singularities.  
\begin{theorem}\label{thm:main3}
If $g$ is the geometric monodromy of an isolated plane curve
singularity, then
$$
\HF(g,+) = 0.
$$
\end{theorem}
The rest of the article is organized as follows.
In Section~\ref{se:monofloer} we recall the basic facts about monotone
symplectomorphisms and Floer homology.  
For background material on symplectic Floer homology in two dimensions 
we refer to \cite{S2} and the references given therein. 
Section~\ref{se:diff} is devoted to diffeomorphisms of finite type and their
properties relevant for Floer homology. We compute the fixed point classes,
establish monotonicity and show that the symplectic action is exact. 
At several points we use a topological proposition about
products of disjoint Dehn twists which was already used in \cite{S1}. We give
a proof of this proposition in Appendix~\ref{ap:dehn}. The result about the
fixed point classes, Proposition~\ref{prop:fclass}, is already contained in
the work of Jiang and Guo~\cite{JG} on the Nielsen number of surface
diffeomorphisms. For the sake of completeness, we give an independent proof.
In Section~\ref{se:main}, the results from the previous sections are put
together to prove Theorem~\ref{thm:main1}.

We would like to point out that the proof of Theorem 1 is very much inspired
by Seidel's computation of the Floer homology of a Dehn twist \cite{S1}. 
In particular our result about the
Floer connecting orbits, Proposition~\ref{prop:corbits}, is a generalization
of the corresponding result for Dehn twists \cite[Lemma 4]{S1}. Our approach
to the proof, however, is slightly 
different from Seidel's original approach. This has led to a shorter proof and
to a wider range of application, see the remark at the end of
Section~\ref{se:corbits}.  
\begin{acknowledgement}
The idea of looking at algebraically finite mapping classes was the fruit of a
week of stimulating discussions with Paul Seidel at Ecole polytechnique. 
I am indebted to him for devoting time and sharing his insight. 
This is a part of my PhD~thesis and I am very grateful to my supervisor Dietmar
Salamon for all his advise during the preparation of this work.  
I thank Eduard Zehnder for many helpful suggestions. 
\end{acknowledgement}

\section{Monotonicity and Floer homology}\label{se:monofloer}

In this section we discuss the notion of monotonicity as defined in \cite{S2}
and give an outline of its significance for Floer homology in two
dimensions. For a more detailed account we refer to the original article.  
At the end of this section we give two criteria for monotonicity which we
use in the next section. 
Throughout this article, $\Sigma$ denotes a closed connected and oriented
2-manifold of genus $\geq2$. In this section, we also fix an area form
$\omega$ on $\Sigma$. 

Let $\phi\in\symp(\Sigma,\omega)$, the group of $\omega$-preserving
diffeomorphisms of $\Sigma$. 
The mapping torus of $\phi$,
$$
\tf = \R\times\Sigma/(t+1,x)\sim(t,\phi(x)),
$$
is a 3-manifold fibered over $S^1=\R/\Z$. 
There are two natural second
cohomology classes on $\tf$, denoted by $[\ophi]$ and $c_\phi$. The first one is
represented by the closed two-form $\ophi$ which is induced from the pullback
of $\omega$ to $\R\times\Sigma$. The second is the Euler class of the
vector bundle 
$$
V_\phi = 
\R\times T\Sigma/(t+1,\xi_x)\sim(t,\de\phi_x\xi_x),
$$
which is of rank 2 and inherits an orientation from $T\Sigma$. 
\begin{definition}[Seidel]
$\phi\in\symp(\Sigma,\omega)$ is called {\bf monotone}, if
$$
[\omega_\phi] = (\area_\omega(\Sigma)/\chi(\Sigma))\cdot c_\phi
$$
in $H^2(\tf;\R)$;
$\symp^m(\Sigma,\omega)$ denotes the set of
monotone symplectomorphisms. 
\end{definition}
Now $H^2(\tf;\R)$ fits into the following short exact sequence
\begin{equation}\label{eq:cohomology}
0 \To  \frac{H^1(\Sigma;\R)}{\im(\id-\phi^*)}
\stackrel{\delta}{\To} H^2(\tf;\R)
\stackrel{\iota^*}{\To} H^2(\Sigma;\R)
\To 0.
\end{equation}
The map $\delta$ is defined as follows.
Let $\rho:\I\to\R$ be a smooth function which vanishes near $0$ and $1$ and
satisfies $\int_0^1\!\rho\,\de t=1$.
If $\theta$ is a closed 1-form on $\Sigma$, then
$\rho\cdot\theta\wedge\de t$ defines a closed 2-form on $T_\phi$; indeed
$$
\delta[\theta] = [\rho\cdot\theta\wedge\de t].
$$
The map $\iota:\Sigma\inclusion T_\phi$ assigns to each $x\in\Sigma$ the
equivalence class of $(1/2,x)$.
Note, that $\iota^*\ophi=\omega$ and $\iota^*c_\phi$ is the Euler class of
$T\Sigma$. 
Hence, by \eqref{eq:cohomology}, there exists a unique class
$m(\phi)\in H^1(\Sigma;\R)/\im(\id-\phi^*)$ satisfying
$$
\delta\,m(\phi) = [\ophi]-(\area_\omega(\Sigma)/\chi(\Sigma))\cdot c_\phi,
$$
where $\chi$ denotes the Euler characteristic. 
Therefore, $\phi$ is monotone if and only if $m(\phi)=0$.
\smallskip\\
We recall the fundamental properties of $\symp^m(\Sigma,\omega)$ from \cite{S2}. 
By $\diff^+(\Sigma)$, we denote the group of orientation preserving
diffeomorphisms of $\Sigma$. 
\smallskip\\
(Naturality)\label{page:natur}
If $\phi\in\symp^m(\Sigma,\omega),\psi\in\diff^+(\Sigma)$, then 
$\psi^{-1}\phi\psi\in\symp^m(\Sigma,\psi^*\omega)$.
\smallskip\\
(Isotopy) 
Let $(\psi_t)_{t\in\I}$ be an isotopy in $\symp(\Sigma,\omega)$, i.e. a smooth
path with $\psi_0=\id$. 
Then 
$$
m(\phi\circ\psi_1)=m(\phi)+[\flux(\psi_t)_{t\in\I}]
$$  
in $H^1(\Sigma;\R)/\im(\id-\phi^*)$; see \cite[Lemma 6]{S2}.
For the definition of the flux homomorphism see \cite{MS}.
\smallskip\\
(Inclusion) 
The inclusion $\symp^m(\Sigma,\omega)\inclusion\diff^+(\Sigma)$ is a homotopy
equivalence.  
This follows from the isotopy property, surjectivity of the flux homomorphism and
Moser's isotopy theorem \cite{Mo}.
Furthermore, the Earl-Eells Theorem \cite{EE} implies that every connected
component of $\symp^m(\Sigma,\omega)$ is contractible.
\smallskip\\
(Floer homology) 
To every $\phi\in\symp^m(\Sigma,\omega)$ symplectic Floer homology
theory assigns a $\Z_2$-graded vector space $\HF(\phi)$ over $\Z_2$, with an
additional multiplicative structure, called the quantum cap product, 
$$
H^*(\Sigma;\Z_2)\otimes\HF(\phi)\To\HF(\phi).
$$
Each $\psi\in\diff^+(\Sigma)$ induces an isomorphism
$\HF(\phi)\cong\HF(\psi^{-1}\phi\psi)$ of $H^*(\Sigma;\Z_2)$-modules. 
\smallskip\\
(Invariance) 
If $\phi,\phi'\in\symp^m(\Sigma,\omega)$ are isotopic, then
$\HF(\phi)$ and $\HF(\phi')$ are naturally isomorphic as
$H^*(\Sigma;\Z_2)$-modules. 
This is proven in \cite[Page 7]{S2}. 
\smallskip\\ 
Now let $g$ be a mapping class of $\Sigma$, i.e. an isotopy class of $\diff^+(\Sigma)$. 
Pick an area form $\omega$ and a
representative $\phi\in\symp^m(\Sigma,\omega)$ of $g$. 
Then $\HF(\phi)$ is an invariant of $g$, which is
denoted by $\HF(g)$. Note that $\HF(g)$ is independent of the choice of an area
form $\omega$ by Moser's isotopy theorem \cite{Mo} and naturality of Floer homology. 

Let $\phi\in\symp^m(\Sigma,\omega)$. 
We give a brief outline of the definition of $\HF(\phi)$ in the special case
where all the fixed points of $\phi$ are non-degenerate. 
This means that for all $y\in\fix(\phi)$, $\det(\id-\de\phi_y)\ne0$. 
In particular, it follows that $\fix(\phi)$ is a finite set and the
$\Z_2$-vector space  
$$
\CF(\phi) := \Z_2^{\fix(\phi)}
$$
admits a $\Z_2$-grading with $(-1)^{\deg y}=\sign(\det(\id-\de\phi_y))$, 
for all $y\in\fix(\phi)$.
The Floer boundary operator is defined as follows. Let $J=(J_t)_{t\in\R}$ 
be a smooth path of $\omega$-compatible complex structures on 
$\Sigma$ such that $J_{t+1}=\phi^*J_t$. For $y^\pm\in\fix(\phi)$, 
let $\M(y^-,y^+;J,\phi)$ denote the space of smooth maps $u:\R^2\to\Sigma$
which satisfy the Floer equations 
\begin{equation}\label{eq:corbit}
\left\{\begin{array}{l}
u(s,t) = \phi(u(s,t+1)), \\
\p_s u + J_t(u)\p_t u = 0, \\
\lim_{s\to\pm\infty}u(s,t) = y^\pm.
\end{array}\right.
\end{equation}
One way to think of the Floer equations is in terms of the symplectic action. 
Let 
$\Of=\{y\in C^{\infty}(\R,\Sigma)\,|\,y(t)=\phi(y(t+1))\}$ denote the twisted
loop space. The action form is the one-form $\alpha_\omega$ on 
$\Of$ defined by 
\begin{equation}\label{eq:daction}
\alpha_\omega(y)\xi = \int_0^1\omega\big(\frac{\de y}{\de t}(t),\xi(t)\big)\,\de t,
\end{equation}
where $y\in\Of$ and $\xi\in T_y\Of$, i.e. $\xi(t)\in T_{y(t)}\Sigma$ and  
$\xi(t)=\de\phi_{y(t+1)}\xi(t+1)$ for all $t\in\R$.
If $\xi,\xi'\in T_y\Of$, then 
$\int_0^1\omega(\xi'(t),J_t\xi(t))\de t$ defines a metric 
on $\Of$. The negative gradient lines of $\alpha_\omega$ with respect to this
metric are solutions of \eqref{eq:corbit}.

Now to every $u\in\M(y^-,y^+;J,\phi)$ is associated a Fredholm operator
$\De_u$ which linearizes (\ref{eq:corbit}) in suitable Sobolev spaces. The
index of this operator is given by the so called Maslov index $\mu(u)$,
which satisfies $\mu(u)=\deg(y^+)-\deg(y^-)\text{ mod }2$. For a generic
$J$, every $u\in\M(y^-,y^+;J,\phi)$ is regular, meaning that $\De_u$ is onto. 
Hence, by the implicit function theorem, $\M_k(y^-,y^+;J,\phi)$ is
a smooth $k$-dimensional manifold, where $\M_k(y^-,y^+;J,\phi)$ denotes the
subset of those $u\in\M(y^-,y^+;J,\phi)$ with $\mu(u)=k\in\Z$.

Translation of the $s$-variable defines a free $\R$-action on
$\M_1(y^-,y^+;J,\phi)$ and hence the quotient is a discrete set of points.
Assume for the moment that for all $y^\pm\in\fix(\phi)$ this quotient is a
finite set and let $n(y^-,y^+)\in\Z_2$ denote its cardinality mod 2.
Define the linear map 
$$
\p_J:\CF(\phi)\To CF_{*+1}(\phi)\quad\text{by}\quad
\fix(\phi)\ni y \Mapsto \sum_{z} n(y,z)z.
$$
That $\p_J$ is of degree $1$ follows from the equation relating the index
and the degree. That $\p_J$ is a boundary operator, i.e. that
$\p_J\circ\p_J=0$, is due to the so-called gluing theorem. 
For this theorem to hold, as well as for $\M_1(y^-,y^+;J,\phi)/\R$ to
be a finite set, one needs certain bounds on the energy. 
Note that bubbling is not an issue here, since $\pi_2(\Sigma)=0$.
The energy of a map
$u:\R^2\to\Sigma$ is given by 
$$
E(u) = \int_{\R}\int_0^1 \omega\big(\p_tu(s,t),J_t\p_tu(s,t)\big)\,\de t\de s.
$$
It is proven in \cite[Lemma 9]{S2} that if $\phi$ is monotone, then the energy
is constant on each $\M_k(y^-,y^+;J,\phi)$. It follows that
$(\CF(\phi),\p_J)$ is a chain complex and that its homology is an 
invariant of $\phi$, denoted by $\HF(\phi)$, i.e. it is independent of $J$. 

Next we introduce the quantum cap product on $\HF(\phi)$. For this, 
choose a Morse function $f:\Sigma\to\R$ and set 
$$
\CoM(f) := \Z_2^{\crit(f)},
$$
with a $\Z$-grading given by the Morse index $\mind_f$. Choose a Riemannian metric on
$\Sigma$ such that $\nabla\!f$ is a Morse-Smale vector field. 
If $x^\pm\in\crit(f)$ and $\mind_f(x^+)=\mind_f(x^-)+1$, denote
by $l(x^-,x^+)\in\Z_2$ the number mod 2 of positive gradient lines going from
$x^-$ to $x^+$. Define the Morse coboundary operator 
$$
\delta_{\nabla\!f}:\CoM(f)\To CM^{*+1}(f)\quad\text{by}\quad
\crit(f)\ni x \Mapsto \sum_y l(x,y)y.
$$
The cohomology of $(\CoM(f),\delta_{\nabla\!f})$ is naturally isomorphic to
$H^*(\Sigma;\Z_2)$, see \cite{Sc}. 
Now by a suitable choice of the function $f$ or the metric, we may assume that for all
$y^\pm\in\fix(\phi),x\in\crit(f)$ and $k\in\Z$, the evaluation map 
$$
\eta_k:\M_k(y^-,y^+;J,\phi)\To\Sigma,\quad
u\Mapsto u(0,0),
$$
is transverse to the unstable manifold $W^u(\nabla\!f,x)\subset\Sigma$. 
Note that the dimension of $W^u(\nabla\!f,x)$ is $2-\mind_f(x)$. 
Hence, if $k=\mind_f(x)$, then $\eta_k^{-1}(W^u(f,x))$ is a discrete set of
points. 
It is in fact a finite set, see \cite[page 8]{S2}. To prove this one 
uses the Gromov-Floer compactification of the moduli spaces and the fact that
$\pi_2(\Sigma)=0$. 
Denote by $q(x;y^-,y^+)\in\Z_2$ the cardinality mod 2 of 
$\eta_k^{-1}(W^u(f,x))\subset\M_k(y^-,y^+;J,\phi)$, where $k=\mind_f(x)$.
Define the linear map
\begin{equation}\label{eq:qca}
\CoM(f)\otimes\CF(\phi)\To\CF(\phi),\quad
x\otimes y\Mapsto\sum_{z}q(x;y,z)z.
\end{equation}
It can be shown that this is a chain map and that the induced map on homology
is independent of $\nabla\!f$ and $J$. It is called the quantum cap product.
For details we refer to \cite{S2} and the references given therein. 

If $\phi$ has degenerate fixed points one needs to perturb equations 
\eqref{eq:corbit} in order to define the Floer homology. Equivalently, one
could say that the action form needs to be perturbed. 
At this point Seidel's approach differs from the usual one. He uses a larger
class of perturbations, but such that the perturbed action form is still
cohomologous to the unperturbed. As a consequence, the usual invariance of
Floer homology under Hamiltonian isotopies is extended to the stronger
property stated above. This ends the general discussion of Floer homology .

To compute the Floer homology of a mapping class one needs to pick a
monotone representative. We now give two criteria for monotonicity which we
use later on. Let $\omega$ be an area form on $\Sigma$ and
$\phi\in\symp(\Sigma,\omega)$. 
\begin{lemma}\label{lemma:monotone1}
Assume that every class $\alpha\in\ker(\id-\phi_*)\subset H_{1}(\Sigma;\Z)$ is
represented by a map $\gamma:S\to\fix(\phi)$, where $S$ is a compact oriented
1-manifold. Then $\phi$ is monotone.
\end{lemma}
\begin{proof}
By dualizing the exact sequence~\eqref{eq:cohomology}, we get the following
exact sequence for homology with real coefficients
\begin{equation}\label{eq:homology}
0\To
H_2(\Sigma;\R)\stackrel{\iota_*}{\To}
H_2(\tf;\R)\stackrel{\hat{\p}}{\To}
\ker(\id-\phi_*)\subset H_1(\Sigma;\R),
\end{equation}
where $\hat{\p}$ is dual to $\delta$. 
Hence, $\phi$ is monotone if and only if 
$$
\langle m(\phi),\alpha\rangle=0, \quad
\forall\;\alpha\in\ker(\id-\phi_*)\subset H_{1}(\Sigma;\R).
$$
If we think of $H_1(\Sigma;\Z)$ as a lattice in $H_1(\Sigma;\R)$, it is
furthermore enough to consider $\alpha\in H_1(\Sigma;\Z)\cap\ker(\id-\phi_*)$.
\smallskip\\
Let $\gamma:S\to\fix(\phi)$ and define $u:S\times S^1\to\tf$ 
by $u(s,t)=(t,\gamma(s))$.
From the definition of $\delta$ on page~\pageref{eq:cohomology}, it is
straight forward to check that 
$$
\langle\delta\alpha,[u]\rangle = \langle\alpha,[\gamma]\rangle,
$$
for all $\alpha\in H_1(\Sigma;\R)$, i.e. that $\hat{\p}[u]=[\gamma]$.
Here, the brackets denote homology classes. 
Now on one hand, since $\p_t u(s,t)=(1,0)$, we have that $u^*\ophi=0$ and hence
$\langle[\ophi],[u]\rangle=0$. 
On the other hand, $\langle c_\phi,[u]\rangle=0$. This is because the
bundle $u^*V_\phi$ is isomorphic to the bundle $\gamma^*T\Sigma\times S^1$,
which is trivial. Hence, it follows that 
$$
\langle m(\phi),[\gamma]\rangle 
= \langle[\ophi],[u]\rangle-
(\area_\omega(\Sigma)/\chi(\Sigma))\langle c_\phi,[u]\rangle 
= 0.
$$
This proves the lemma. 
\end{proof}
\begin{lemma}\label{lemma:monotone2}
If $\phi^k$ is monotone for some $k>0$, then $\phi$ is monotone. 
If $\phi$ is monotone, then $\phi^k$ is monotone for all $k>0$. 
\end{lemma}
\begin{proof}
Recall that $\tf$ is the orbit space of the $\Z$-action
$n\cdot(t,x) = (t+n,\phi^{-n}(x))$, where $n\in\Z$ and $(t,x)\in\R\times\Sigma$.
If we only divide out by the 
subgroup $k\Z$, for $k\in\N_{>0}$, we naturally get the mapping torus of
$\phi^k$. 
Further dividing by $\Z/k\Z$ defines the $k$-fold covering map
$p_k:T_{\phi^k}\to\tf$. 
It is straight forward to check that 
\begin{equation}
p_k^*[\omega_\phi] = [\omega_{\phi^k}] \quad\text{and}\quad 
p_k^*c_\phi = c_{\phi^k}.
\label{eq:iterate}\end{equation}
The first equality follows immediately from the definitions. 
To prove the second, note that
$$  
p_k^*\big((TM\times\R)/\Z\big) \cong (TM\times\R)/k\Z \cong V_{\phi^k}, 
$$
where the $\Z$-action on $\R\times T\Sigma$ is given by 
$n\cdot(t,\xi_x)=(t+n,\de\phi_x^{-n}\xi_x)$, for $n\in\Z$ and $\xi_x\in
T_x\Sigma$. 
The lemma follows from \eqref{eq:iterate} and the fact that $p_k^*$ is
injective. To prove injectivity, define the map
$a^k_*:H^2(T_{\phi^k};\R)\to H^2(T_\phi;\R)$ by averaging differential forms;
$a^k_*$ is a left inverse of $p_k^*$, i.e. $a^k_*\circ p_k^*=\id$. This ends
the proof of the lemma.
\end{proof}

\section{Diffeomorphisms of finite type}\label{se:diff}

We begin with the basic definition. Note that $S^1$ is always identified 
with $\R/\Z$.
\begin{definition}\label{def:ftype}
We call $\phi\in\diff_+(\Sigma)$ of {\bf finite type} if the following holds.
There is a $\phi$-invariant finite union $N\subset\Sigma$ of disjoint
non-contractible annuli such that:
\smallskip\\
(1) $\phi|\Sigma\setminus N$ is periodic, i.e. there exists
$\ell>0$ such that $\phi^\ell|\Sigma\setminus N=\id$.
\smallskip\\
(2) Let $N'$ be a connected component of $N$ and $\ell'>0$ be the
smallest integer such that $\phi^{\ell'}$ maps $N'$ to itself. Then
$\phi^{\ell'}|N'$ is given by one of the following two models with respect to
some coordinates $(q,p)\in\I\times S^1$:
\medskip\\
\begin{minipage}{4cm}
(twist map)
\end{minipage}
\begin{minipage}{6cm}
$(q,p)\Mapsto(q,p-f(q))$ 
\end{minipage}
\medskip\\
\begin{minipage}{4cm}
(flip-twist map)
\end{minipage}
\begin{minipage}{6cm}
$(q,p)\Mapsto(1-q,-p-f(q))$,
\end{minipage}
\medskip\\
where $f:\I\to\R$ is smooth and strictly monotone. 
A twist map is called {\bf positive/negative}, 
if $f$ is increasing/decreasing. 
\smallskip\\ 
(3) Let $N'$ and $\ell'$ be as in (2). 
If $\ell'=1$ and $\phi|N'$ is a twist map, then $\im(f)\subset[0,1]$,
i.e. $\phi|\text{int}(N')$ has no fixed points.  
\smallskip\\
(4) If two connected components of $N$ are homotopic, then the corresponding local
models of $\phi$ are either both positive or both negative twists.
\end{definition}
\begin{remarks}
(i) 
Let $\phi$ be a diffeomorphism of finite type and 
$\ell$ be as in (1). 
Then $\phi^\ell$ is the product of (multiple) {\bf Dehn twists} ``along $N$''.
Moreover, two parallel Dehn twists have the same sign, by (4). We say that
$\phi$ has {\bf uniform twists}, if $\phi^\ell$ is the product of only
positive, or only negative Dehn twists. 
\smallskip\\ 
(ii) A mapping class of $\Sigma$ is called {\bf algebraically finite} if it
does not have any pseudo-Anosov components in the sense of Thurston's
theory of surface diffeomorphism. Every such class is represented by a
diffeomorphism of finite type. 

To see this, recall Thurston's classification
theorem, \cite[Theorem 4]{Th}: 
for every mapping class of $\Sigma$, there exists a diffeomorphism $\phi$
representing the class and a $\phi$-invariant finite union $C\subset\Sigma$ of
non-contractible disjoint circles such that: 
\smallskip\\
$(1')$ The components of $C$ are pairwise non-homotopic, 
\smallskip\\
$(2')$ If $\Sigma'$ is a $\phi$-invariant union of connected components of 
$\Sigma\setminus C$, then $\phi|\Sigma'$ is isotopic to either a periodic
or a pseudo-Anosov map.    
\smallskip\\
The set $C$ is called a reducing set. 
Starting with a mapping class without pseudo-Anosov components, one
first chooses a minimal reducing set $C$, meaning that it has
the minimal number of components of all reducing sets. 
Minimality guarantees that after isotopying the Nielsen-Thurston
representative $\phi$ on a complement of a tubular neighborhood $N$ of $C$ to
a periodic map, $\phi|N$ does not have periodic components. One can thus
achieve condition (2) above, by isotopying $\phi|N$ relative to $\p N$.
If (3) is not satisfied, this is achieved in a last step by introducing
further components of $C$, violating $(1')$, but such that (4) still
holds. 
\smallskip\\ 
(iii)
The term algebraically finite goes back to Nielsen~\cite{Ni}. 
Fried \cite{Fr} defined the notion of algebraically finite
diffeomorphism in any dimension. In two dimensions, these are
special representatives of algebraically finite mapping classes.
Fried's definition, however, is adopted to the theory of dynamical systems. 
For our purpose, a representative which is of the special type defined above
is most convenient. 
\smallskip\\ 
(iv) 
The term flip-twist map is taken from \cite{JG}. 
\end{remarks}
The rest of this section is devoted to the study of diffeomorphisms of finite
type. The points of interest for Floer homology are: fixed point classes,
monotonicity and action. The results we obtain are used in
Section~\ref{se:main} to compute the Floer homology.  
\begin{convention} From now on, $\phi$ denotes a diffeomorphism of
finite type and $N$ the associated $\phi$-invariant union of annuli. 
By $\Sigma_0$ we denote the union of the components
of $\Sigma\setminus\text{int}(N)$, where $\phi$ restricts to the identity.
Furthermore, we denote by $\ell$ the smallest positive integer such that
$\phi^\ell$ restricts to the identity on $\Sigma\setminus N$. 
\end{convention}
The first proposition describes the set of fixed point classes of $\phi$. It
is a special case of a theorem by B.~Jiang and J.~Guo~\cite{JG}, which gives
for any mapping class a representative that realizes its
Nielsen number. 
\begin{prop}[Fixed point classes]\label{prop:fclass}
Each fixed point class of $\phi$ is either a connected component of $\Sigma_0$ or
consists of a single fixed point. A fixed point $x$ of the second type satisfies
$\det(\id-\de\phi_x)>0$.  
\end{prop}
The crucial step in our proof of this proposition is to prove it in the
special case of products of disjoint Dehn twists. For this, we refer to
Appendix~\ref{ap:dehn}.  
\begin{proof}[Proof of Proposition \ref{prop:fclass}]
First note, that if $x\in\fix(\phi)\cap\text{int}(N)$, then $\phi$
restricted to the component of $N$ containing $x$, is a flip-twist map and
$x=(\frac{1}{2},-\frac{1}{2}f(\frac{1}{2}))$ or
$(\frac{1}{2},\frac{1}{2}-\frac{1}{2}f(\frac{1}{2}))$.  
Now let $x\ne y$ be arbitrary fixed points in the same fixed point class. 
We prove in three steps, that $x$ and $y$ are in the same connected component
of $\Sigma_0$. 
\begin{claim}{1} 
$x$ and $y$ are in the same component of either $\text{int}(N)$ or
$\Sigma\setminus\text{int}(N)$.
\end{claim}
\begin{pc}
Note that every connected component of $\Sigma\setminus\text{int}(N)$ is a
connected component of $\fix(\phi^\ell)$. Similarly, if $x\in\text{int}(N)$, 
then $x$ is contained in a fixed circle of $\phi^\ell$. 
Such a circle is also a connected component of $\fix(\phi^\ell)$.
By Corollary~\ref{cor:dehn} in Appendix~\ref{ap:dehn}, however, 
every connected component of $\fix(\phi^\ell)$ is a fixed point class of
$\phi^\ell$. 
Since $x$ and $y$ are in the same fixed point class of
$\phi^\ell$, this proves claim 1. Denote by $M$ the connected component of $N$ or
$\Sigma\setminus\text{int}(N)$, containing $x$ and $y$.
\end{pc}
\begin{claim}{2}
$x$ and $y$ are in the same fixed point class of $\phi|M$. 
\end{claim}
\begin{pc}
By assumption, there exists a map $u:\I^2\to\Sigma$ with 
$$
\quad u(0,t) = x,\quad u(1,t) = y \quad\text{and}\quad
u(s,1)=\phi(u(s,0)),
$$
for all $s,t\in[0,1]$.
By Corollary~\ref{cor:dehn}, we can assume that 
$u(s,0)\in M$, for all $s\in\I$. 
Let $M'$ be the union of $M$ and a tubular neighborhood of $\p M$. 
We prove that $u$ can be deformed in the
interior of $\I^2$ such that its image is contained in a $M'$. 
First note, that by a small
perturbation, we may assume that $u$ is transverse to $\p M'$. 
Hence $u^{-1}(\p M')\subset\I^2$ is a 1-dimensional
submanifold with boundary $\p u^{-1}(\p M')=u^{-1}(\p M')\cap\p\I^2=\emptyset$. 
Every component of $u^{-1}(\p M')$ is therefore a circle and bounds a disk
in $\I^2$. The restriction of $u$ to this such a disk represents an element of
$\pi_2(\Sigma,\p M')$. Since $\pi_2(\Sigma,\p M')=0$, $u$ can be deformed
in the interior of $\I^2$ to a map $v$ such that the number of components of
$v^{-1}(\p M')$ is less than that of $u^{-1}(\p M')$. 
It follows inductively, that $u$ can be deformed in the interior
of $\I^2$ to a map $w$ with $w^{-1}(\p M')=\emptyset$. This proves
claim 2 and we are left with
\end{pc}
\begin{claim}{3}
Let $\varphi$ be either a flip-twist map or a non-trivial orientation
preserving periodic diffeomorphism of a compact connected surface $M$
of Euler characteristic $\leq0$. Then each fixed point class of $\varphi$
consists of a single point. 
\end{claim}
\begin{pc}
For a flip-twist map, this is checked explicitly by using the model. 
The other case was first proven in \cite{J}. We repeat the
argument here.
First assume that $M$ is closed. The uniformization theorem states 
that in every conformal class of metrics on $M$, there is a unique metric of
constant curvature $-1$ if $\chi(M)<0$ or 0 if $\chi(M)=0$. This implies that
the unique representative of a $\varphi$-invariant conformal class, such a
class exists since $\varphi$ is finite order, is itself $\varphi$-invariant. Hence
we can pick a $\varphi$-invariant metric of constant curvature $-1$ or 0 on $M$ and 
lift $\varphi$ to an isometry $\tilde{\varphi}$ of the universal cover
$\tilde{M}$ of $M$. $\tilde{M}$ is either isometric to the hyperbolic plane
$\Hyp^2$ or the Euclidean plane $\R^2$. 

Let $x\in\fix(\varphi)$ and 
let $\tilde{\varphi},\tilde{x}$ be lifts of $\varphi,x$ to
$\tilde{M}$, such that $\tilde{\varphi}(\tilde{x})=\tilde{x}$. 
Note, that a fixed point of $\varphi$ is in the same class as $x$ if and only
if it can be lifted 
to a fixed point of $\tilde{\varphi}$. Assume by contradiction that
$\tilde{y}\ne\tilde{x}$ is a fixed point of $\tilde{\varphi}$. It follows that
the unique geodesic going through $\tilde{x}$ and $\tilde{y}$ is pointwise fixed by
$\tilde{\varphi}$. In particular, since $\tilde{\varphi}$ preserves
orientation, $\de\tilde{\varphi}_{\tilde{x}}=\id$. This implies
that $\tilde{\varphi}=\id$, because an isometry of $\Hyp^2$ or $\R^2$ is
determined by its value and differential at one point. This proves claim 3 in
the case that $M$ is closed.  

The case $\p M\ne\emptyset$ is reduced to the above case by gluing two
copies of $M$ together along a $\varphi$-invariant tubular neighborhood of $\p
M$. The glued manifold is closed and of Euler characteristic $\leq0$;
$\varphi$ extends to a non-trivial diffeomorphism $\varphi'$, which is
orientation preserving and of finite order.
Hence, every fixed point class of $\varphi'$ is a single point.
The same therefore holds for $\varphi$. This ends the proof of claim 3. 
Finally we have 
\end{pc}
\begin{claim}{4}
If $x\in\fix(\phi)\setminus\Sigma_0$, then $\det(\id-\de\phi_x)>0$.
\end{claim}
\begin{pc}
The point $x$ is a fixed point of either a flip-twist map or an orientation
preserving non-trivial isometry. 
In the first case, the assertion is checked by using the local model. 
Similarly, one checks in the second case, that 
$\det(\id-\de\phi_x)\leq0$ if and only if $\de\phi_x=\id$. As shown in the
proof of claim 3, however, $\de\phi_x=\id$ implies that $x\in\Sigma_0$, which
is a contradiction.  
\end{pc}
\end{proof}
The next issue is monotonicity. First note that
if $\omega'$ is an area form on $\Sigma$ which is the standard form
$\de q\wedge\de p$ with respect to the $(q,p)$-coordinates on $N$, then 
$\omega:=\sum_{i=1}^\ell(\phi^i)^*\omega'$ is 
standard on $N$ and $\phi$-invariant, i.e. $\phi\in\symp(\Sigma,\omega)$.
To prove that $\omega$ can be chosen such that $\phi\in\symp^m(\Sigma,\omega)$,
we distinguish two cases: uniform and non-uniform twists. In the first case we
have the following stronger statement. 
\begin{prop}[Monotone 1]\label{prop:monotone3}
If $\phi$ has uniform twists and $\omega$ is a $\phi$-invariant area
form, then $\phi\in\symp^m(\Sigma,\omega)$.
\end{prop}
\begin{proof}
By Lemma~\ref{lemma:monotone2}, it is enough to prove that $\phi^\ell$ is
monotone with respect to any $\phi$-invariant area form. Replace $\phi$ by
$\phi^\ell$. By the uniform twist condition, $\phi$ is the product of 
disjoint Dehn twists which are all, say positive. 
We prove that $\phi$ satisfies the hypothesis of Lemma~\ref{lemma:monotone1}
and is therefore monotone. 
We use the Picard-Lefschetz formula for the action of a positive Dehn twist on
$H_1(\Sigma;\Z)$: if the twist is along $C\subset\Sigma$, then
$\alpha\mapsto\alpha-(\alpha\cdot[C])[C]$, where $\alpha\in H_1(\Sigma;\Z)$
and $[C]$ denotes the homology class of $C$ with respect to some orientation. 
The dot stands for the intersection pairing. 

Let $C_1,\dots,C_n\subset\Sigma$ be the disjoint non-contractible circles 
along which $\phi$ twists. Choose orientations of the $C_i$. 
Let $\alpha\in\ker(\id-\phi_*)\subset H_1(\Sigma;\Z)$. 
We claim that for all $i=1,\dots,n$, $\alpha\cdot[C_i]=0$.
This is equivalent to the condition that $\alpha$ is represented
by a map $S\to\fix(\phi)$, where $S$ is a compact oriented 1-manifold, and
therefore ends the proof of the proposition.   
Since the $C_i$ are pairwise disjoint, it follows from the Picard-Lefschetz
formula that  
$$ 
\alpha = \phi_*\alpha = \alpha - \sum_{i=1}^n (\alpha\cdot[C_i])[C_i]
$$
and hence that $\sum_{i=1}^n (\alpha\cdot[C_i])[C_i]=0$. 
Pairing with $\alpha$, we get
$\sum_{i=1}^n (\alpha\cdot[C_i])^2 = 0$,
which implies that $\alpha\cdot[C_i]=0$, for all $i=1,\dots,n$.
\end{proof}
In the non-uniform case, monotonicity is a more subtle point and does not 
hold for arbitrary $\phi$-invariant area forms. 
\begin{prop}[Monotone 2]\label{prop:monotone4}
If $\phi$ does not have uniform twists, there exists a $\phi$-invariant
area form $\omega$ such that $\phi\in\symp^m(\Sigma,\omega)$. Moreover,
$\omega$ can be chosen such that it is the standard form $\de q\wedge\de p$ on
$N$.  
\end{prop}
\begin{proof}
The strategy of the proof is the following.
Assume first that $\phi|\Sigma\setminus N=\id$. 
We begin by defining a $\phi$-invariant area form $\omega$ with
$\area_\omega(\Sigma)=-\chi(\Sigma)$.
Then we construct for every class $[\gamma]\in\ker(\id-\phi_*)$, a class
$[\Gamma]\in H_2(\tf;\Z)$ such that 
$$
\hat{\p}[\Gamma]=[\gamma]
\quad\text{and}\quad 
\langle[\ophi],[\Gamma]\rangle=-\langle c_\phi,[\Gamma]\rangle,
$$ 
with $\hat{\p}$ as defined in the sequence~\eqref{eq:homology}. 
Finally, we show how the general case is reduced to the case above.

We start with the following set-up. 
Fix a union $\tilde{N}\subset\Sigma$ of
disjoint non-contractible and pairwise non-homotopic annuli such that
$\phi|\Sigma\setminus\Nn=\id$.
Moreover, let for every 
connected component $N'$ of $\tilde{N}$, $\ell'$ be a positive integer and
$f:[0,1]\to\R$ be a smooth monotone function 
with $f(0)=0,f(1)=\ell'$ and such that $\phi|N'$ is
(an $\ell'$-fold Dehn twist) given by the model 
$(q,p)\mapsto(q,p-f(q))$ for $(q,p)\in\I\times S^1$. 
We emphasize here, that not only the function $f$ but also the local chart of
$N'$ is fixed for the rest of the proof.

In a first step, we choose for every component $N'$ of $\Nn$, an embedded circle
$C'\subset N'$ as follows.
Look at the set $\text{graph}(-f)\subset\I\!\times\![-\ell',0]$ if $\ell'>0$,
respectively $\I\!\times\![0,-\ell']$ if $\ell'<0$. See the figures below.
\bigskip\\
\begin{minipage}{6cm}
\begin{center}
\begin{picture}(0,0)%
\includegraphics{monotone2.pstex}%
\end{picture}%
\setlength{\unitlength}{987sp}%
\begingroup\makeatletter\ifx\SetFigFont\undefined%
\gdef\SetFigFont#1#2#3#4#5{%
  \reset@font\fontsize{#1}{#2pt}%
  \fontfamily{#3}\fontseries{#4}\fontshape{#5}%
  \selectfont}%
\fi\endgroup%
\begin{picture}(5100,5643)(3976,-6373)
\put(5776,-2161){\makebox(0,0)[lb]{\smash{\SetFigFont{5}{6.0}{\familydefault}{\mddefault}{\updefault}{\color[rgb]{0,0,0}{\small $-$}}%
}}}
\put(6751,-5836){\makebox(0,0)[lb]{\smash{\SetFigFont{5}{6.0}{\familydefault}{\mddefault}{\updefault}{\color[rgb]{0,0,0}{\small $+$}}%
}}}
\put(6451,-1336){\makebox(0,0)[lb]{\smash{\SetFigFont{5}{6.0}{\rmdefault}{\mddefault}{\updefault}{\color[rgb]{0,0,0}{\small $a$}}%
}}}
\put(8326,-1336){\makebox(0,0)[lb]{\smash{\SetFigFont{5}{6.0}{\familydefault}{\mddefault}{\updefault}{\color[rgb]{0,0,0}{\small 1}}%
}}}
\put(4501,-886){\makebox(0,0)[lb]{\smash{\SetFigFont{5}{6.0}{\familydefault}{\mddefault}{\updefault}{\color[rgb]{0,0,0}{\small $p$}}%
}}}
\put(9076,-1411){\makebox(0,0)[lb]{\smash{\SetFigFont{5}{6.0}{\familydefault}{\mddefault}{\updefault}{\color[rgb]{0,0,0}{\small $q$}}%
}}}
\put(3976,-6286){\makebox(0,0)[lb]{\smash{\SetFigFont{5}{6.0}{\familydefault}{\mddefault}{\updefault}{\color[rgb]{0,0,0}{\small $-\ell'$}}%
}}}
\end{picture}
\end{center}
\par\medskip
\begin{minipage}{6cm}
\begin{center} 
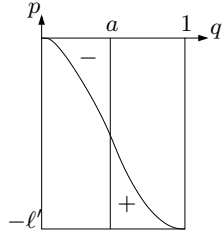
Figure 1: $\text{graph}(-f)$ \\
for a positive twist
\end{center}
\end{minipage}
\end{minipage}
\hfill
\begin{minipage}{6cm}
\begin{center}
\begin{picture}(0,0)%
\includegraphics{monotone1.pstex}%
\end{picture}%
\setlength{\unitlength}{987sp}%
\begingroup\makeatletter\ifx\SetFigFont\undefined%
\gdef\SetFigFont#1#2#3#4#5{%
  \reset@font\fontsize{#1}{#2pt}%
  \fontfamily{#3}\fontseries{#4}\fontshape{#5}%
  \selectfont}%
\fi\endgroup%
\begin{picture}(5250,6081)(3826,-6811)
\put(7021,-2191){\makebox(0,0)[lb]{\smash{\SetFigFont{5}{6.0}{\familydefault}{\mddefault}{\updefault}{\color[rgb]{0,0,0}{\small $-$}}%
}}}
\put(6061,-5851){\makebox(0,0)[lb]{\smash{\SetFigFont{5}{6.0}{\familydefault}{\mddefault}{\updefault}{\color[rgb]{0,0,0}{\small $+$}}%
}}}
\put(6556,-6736){\makebox(0,0)[lb]{\smash{\SetFigFont{5}{6.0}{\familydefault}{\mddefault}{\updefault}{\color[rgb]{0,0,0}{\small $a$}}%
}}}
\put(4501,-886){\makebox(0,0)[lb]{\smash{\SetFigFont{5}{6.0}{\familydefault}{\mddefault}{\updefault}{\color[rgb]{0,0,0}{\small $p$}}%
}}}
\put(9076,-6211){\makebox(0,0)[lb]{\smash{\SetFigFont{5}{6.0}{\familydefault}{\mddefault}{\updefault}{\color[rgb]{0,0,0}{\small $q$}}%
}}}
\put(8326,-6811){\makebox(0,0)[lb]{\smash{\SetFigFont{5}{6.0}{\familydefault}{\mddefault}{\updefault}{\color[rgb]{0,0,0}{\small 1}}%
}}}
\put(3826,-1786){\makebox(0,0)[lb]{\smash{\SetFigFont{5}{6.0}{\familydefault}{\mddefault}{\updefault}{\color[rgb]{0,0,0}{\small $-\ell'$}}%
}}}
\end{picture}
\end{center}
\par\medskip
\begin{minipage}{6cm}
\begin{center}
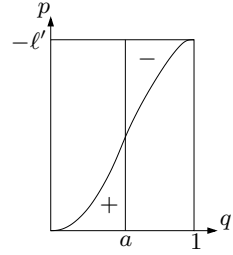
Figure 2: $\text{graph}(-f)$ \\
for a negative twist
\end{center}
\end{minipage}
\end{minipage}
\bigskip\\
For any $a\in(0,1)$, the complement of the union of $\text{graph}(-f)$ and the set
$\{q=a\}$ has four components.  
If $\phi|N'$ is a positive twist, we choose $a$ such that
the left upper component (indicated with a $-$ sign in Figure~1) and the right 
lower component (indicated with a $+$ sign) have the same area with respect to
the standard area form on $\I\!\times\![-\ell',0]$. 
If $\phi|N'$ is a negative twist, left upper is replaced by left
lower and right lower by right upper and the sings are interchanged.
In both cases, we set $C':=\{a\}\times S^1\subset N'$, with orientation
induced from $S^1$.
The purpose of this construction will become clear below. Let $C$ denote the
union of the loops $C'$. Let $\Sigma_1,\dots,\Sigma_m\subset\Sigma$ denote the
closures of the connected components of $\Sigma\setminus C$.
Since the $C'$ are disjoint non-contractible and pairwise non-homotopic, it 
follows that $\chi(\Sigma_j)<0$ for all $j=1,\dots,m$. Now choose an area form
$\omega$ on $\Sigma$ such that 
\begin{equation}\label{eq:mono4area}\begin{split}
\area_\omega(\Sigma_j) 
&= -\chi(\Sigma_j) \quad\text{for all }j=1,\dots,m,\\
\omega|\Nn
&= \e\cdot\de q\wedge\de p,
\end{split}\end{equation}
where $\e>0$ is sufficiently small.
By the first condition, we have that $\area_\omega(\Sigma)=-\chi(\Sigma)$
and from the second it follows that $\phi^*\omega=\omega$.
We now prove in several steps that $[\omega_\phi]=-c_\phi$ in $H^2(T_\phi;\R)$. 

Let $S$ be a compact oriented 1-manifold and $\gamma:S\to\Sigma$ be an
immersion which is transverse to $C$. Moreover, assume that 
$[\gamma]=[\phi\circ\gamma]$ in $H_1(\Sigma;\Z)$.
The goal is to lift the 1-cycle $\gamma$ to a
2-cycle $\Gamma$ in $\tf$. For this, we first define a 2-chain $A$ in $\Sigma$
which satisfies
$$ 
\p A=\gamma-\phi\!\circ\!\gamma-\sum_{i=1}^n 
\ell_i\big([\gamma]\cdot[C_i]\big)C_i,
$$
where we think of the right hand side as a 1-chain. 
Here, we have introduced a numbering of the components of $C$. 
The chain $A$ can be described a follows, compare Figure~1, 2: 
At every intersection point of $\gamma$ and $C$ where $\gamma$ runs in 
the positive $q$-direction, there is a local contribution to
$A$ given by the regions in Figure~1, 2 which are labelled by $\pm$. 
The sign of the contribution is as indicated in the figure. 
If $\gamma$ runs in the negative $q$-direction, the signs are
interchanged. Note that by our choice of
$\omega$, we have that $\int_A\omega=0$. 

Next, we use that $\gamma$ is homologous to $\phi\circ\gamma$, i.e. that
$$
\sum_{i=1}^n \ell_i([\gamma]\!\cdot\![C_i])[C_i]=0. 
$$
This means that there exist
integers $k_1,\dots,k_m$ such that
$$
\p\big(\sum_{j=1}^m k_j\Sigma_j\big) =
\sum_{i=1}^n \ell_i\big([\gamma]\cdot[C_i]\big)C_i.
$$
We can now define the 2-chain $\Gamma$ in
$\tf=\I\times\Sigma/(0,\phi(x))\sim(1,x)$, by
$$
\Gamma:=-[0,1/2]\times(\phi\circ\gamma)
-\{1/2\}\times\big(A+{\textstyle\sum_{j=1}^m k_j\Sigma_j}\big)
-[1/2,1]\times\gamma.
$$
By construction, $\Gamma$ is a cycle; indeed 
\begin{equation*}
\begin{split}
\p\Gamma 
&=
\{0\}\times(\phi\circ\gamma) - \{1/2\}\times(\phi\circ\gamma) 
- \{1/2\}\times\p A \\
&\quad 
{}- \{1/2\}\times\p(\textstyle{\sum_{j=1}^m}k_j\Sigma_j)
+\{1/2\}\times\gamma - \{1\}\times\gamma \\
&= 
0.
\end{split}
\end{equation*}
By a similar calculation as in the proof of
Lemma~\ref{lemma:monotone1}, it furthermore follows that
$\hat{\p}[\Gamma]=[\gamma]$.  
\begin{claim}{1} 
$\langle[\ophi],[\Gamma]\rangle=-\sum_{j=1}^m k_j\area_\omega(\Sigma_j)$.
\end{claim}
\begin{pc}
Only the middle summand of $\Gamma$ contributes to
$\langle[\ophi],[\Gamma]\rangle$. Since $A$ has vanishing
$\omega$-area, this already proves claim 1.
\end{pc}
\begin{claim}{2}
$\langle c_\phi,[\Gamma]\rangle=-\sum_{j=1}^m k_j\chi(\Sigma_j)$.
\end{claim}
\begin{pc}
To prove this, we use the following property of the Euler class. 
If a smooth section $s:T_\phi\to V_\phi$ is transverse to the zero-section, 
then $s^{-1}(0)\subset T_\phi$ is a submanifold of codimension 2 and
its homology class, with respect to a suitable orientation, is
Poincar{\'e}-dual to $c_\phi$.
In particular,
$\langle c_\phi,[u]\rangle$ equals the intersection number
$[s^{-1}(0)]\cdot[u]$, for any $[u]\in H_2(T_{\phi};\Z)$. 
The orientation of $s^{-1}(0)$ at a point $x$ is defined as follows. 
Let $\{e_1,e_2,e_3\}$ be an oriented basis of $T_xT_\phi$ such that $e_1$
is tangent to $s^{-1}(0)$. 
Then $e_1$ is said to be oriented, if $\{e_1,e_2,e_3,\de s_xe_2,\de s_xe_3\}$
is an oriented basis of
$$
T_{(x,0)}V_\phi\cong \R\oplus T_x\Sigma\oplus T_x\Sigma.
$$
We now define a smooth section of $V_\phi$.
To begin with, we choose a vector field $\xi$ on $\Sigma$ with only
non-degenerate zeros and such that $\xi|\Nn=\p/\p q$.
Furthermore, we require that $\xi^{-1}(0)$ is disjoint from $\im(\gamma)$.
That such a vector field exists is a standard result in differential topology.
By the Poincar{\'e}-Hopf Theorem, the sum of
indices of $\xi$ over all zeros in $\Sigma_j$ equals
$\chi(\Sigma_j)$, for all $j=1,\dots,m$. 
\smallskip\\
Note, that the vector field $\phi^*\xi-\xi$ is supported in $\Nn$,
where it is given by $(0,-f'(q))$ with respect to the local model. 
Hence, there exists a smooth path $(\xi_t)_{t\in\R}$ of vector fields
such that $\xi_{t+1}=\phi^*\xi_t$, $\xi_t=\xi$ on $\Sigma\setminus\Nn$ 
and $\xi_t^{-1}(0)=\xi^{-1}(0)$. 
Let the section $s:T_\phi\to V_\phi$ be defined by
$s([t,x]):=[t,\xi_t(x)]$; recall that 
$$
V_\phi = 
\R\times T\Sigma/(t+1,\xi_x)\sim(t,\de\phi_x\xi_x).
$$
By our choice of the vector field $\xi$, $s$ in transverse to the zero-section
and thus, $[s^{-1}(0)]$ is Poincar{\'e}-dual to $c_\phi$. 
Moreover, 
\begin{eqnarray*}
[s^{-1}(0)]\cdot[\Gamma] &=& 
-\sum_{j=1}^m k_j \big(s^{-1}(0)\cdot\Sigma_j\big) \\
&=& -\sum_{j=1}^m k_j\cdot\chi(\Sigma_j).
\end{eqnarray*}
The numbers $s^{-1}(0)\cdot\Sigma_j$ are well defined because $s^{-1}(0)$
intersects $\Sigma_j$ transversally and in the interior of $\Sigma_j$.
Note, that the sign of an intersection point equals the index of $\xi$ at that
point. This proves claim~2. 
\end{pc}
From claim~1, 2 and the first equation in \eqref{eq:mono4area}, 
we conclude that 
$$
\langle[\ophi],[\Gamma]\rangle=-\langle c_\phi,[\Gamma]\rangle,
$$ 
and hence that 
$\phi\in\symp^m(\Sigma,\omega)$. We end the proof of the proposition with the
following observation. 
Let $\phi$ be a diffeomorphism of finite type;
apply the above construction to $\phi^\ell$ and let $\omega$ be an area form
which satisfies \eqref{eq:mono4area}. It follows that
$$
\omega':=\frac{1}{\ell}\sum_{i=1}^\ell(\phi^i)^*\omega
$$ 
also satisfies \eqref{eq:mono4area} and hence that 
$\phi^\ell\in\symp^m(\Sigma,\omega')$. On the other hand
$\phi^*\omega'=\omega'$, 
and therefore $\phi\in\symp^m(\Sigma,\omega')$, by
Lemma~\ref{lemma:monotone2}. 

The idea of the proof is the following. As above, we first replace $\phi$ by
$\phi^\ell$. Then we construct an $\phi$-invariant area form $\omega$ with
$\area_\omega(\Sigma)=-\chi(\Sigma)$. The main step is the following:
for every integer class $[\gamma]\in\ker(\id-\phi_*)$ we construct a class
$[\Gamma]\in H_2(\tf;\Z)$ such that 
$$
\p[\Gamma]=[\gamma]
\quad\text{and}\quad 
\langle[\ophi],[\Gamma]\rangle=-\langle c_\phi,[\Gamma]\rangle.
$$ 
This proves the proposition for $\phi$. Finally, we show that $\omega$ can be
chosen such that it is invariant under the original $\phi$.
By Lemma~\ref{lemma:monotone2} this proves the proposition.
\end{proof}
Next we consider the symplectic action $\alpha_\omega$ on the twisted loop
space $\Of$ and prove that it is exact. See \eqref{eq:daction} for the
definition of $\alpha_\omega$. This result is crucial for the computation of the Floer
homology, in particular for the use of the connecting orbits proposition
proved in the next section. 

We need the following lemma, which holds for any $\phi\in\symp(\Sigma,\omega)$.
First note that a loop in $\Of$ is represented by a map $u:S^1\times\I\to\Sigma$
with 
$$
u(s,t)=\phi(u(s,t+1))
\quad\text{for all }(s,t)\in S^1\times\I.
$$ 
By $[u]$ we denote the homology class of the loop $u$ in $\Of$. 
\begin{lemma}\label{lemma:mu}
Let $u$ and $v$ be two loops in $\Of$. If $u(.,0)$ and $v(.,0)$ are freely
homotopic loops in $\Sigma$, then 
$\langle[\alpha_\omega],[u]\rangle=\langle[\alpha_\omega],[v]\rangle$.
\end{lemma} 
\begin{proof}
Let $w:S^1\times\I\to\Sigma$ be such that $w(.,0)=u(.,0)$ and
$w(.,1)=v(.,0)$. Define the map 
$u':=w^{-1}\#u\#(\phi^{-1}\circ w):S^1\times\I\to\Sigma$ by  
$$
u'(s,t)=\begin{cases}
w(s,1-3t) & \text{if $t\in[0,1/3]$} \\
u(s,3t-1) & \text{if $t\in[1/3,2/3]$} \\
\phi^{-1}(w(s,3t-2)) & \text{if $t\in[2/3,1]$}.
\end{cases}
$$
Since for all $t\in[0,1/3],s\in S^1$,
$u'(s,1-t)=\phi^{-1}(w(s,1-3t))=\phi^{-1}(u'(s,t))$, 
it follows that $u'$ and $u$ are homotopic loops in $\Of$; in particular
$[u]=[u']$.  
Note that since $u'(s,0)=v(s,0)$ and $u'(s,1)=v(s,1)$ for all $s\in S^1$, 
the map $u'\#v^{-1}$ descends to a map $S^1\times S^1\to\Sigma$. 
Now if $w$ is an arbitrary loop in $\Of$, then
$$
\langle[\alpha_\omega],[w]\rangle=-\int_{S^1\times\I}w^*\omega.
$$
Therefore
$$
\langle[\alpha_\omega],[u']\rangle-\langle[\alpha_\omega],[v]\rangle
= -\int_{S^1\times S^1}(u'\#v^{-1})^*\omega = 0. 
$$
In the last equality we use the fact 
that the mapping degree of $u'\#v^{-1}$
vanishes, since the genus of $\Sigma$ is $\geq2$.
This proves the proposition.
\end{proof}
We return to the situation where $\phi$ is a diffeomorphism of finite type. 
The proof of the following proposition relies on the results discussed in
Appendix~\ref{ap:dehn}. 
\begin{prop}[Action]\label{prop:action}
If $\omega$ is a $\phi$-invariant area form, then $\alpha_\omega$ has
vanishing periods.
\end{prop}
\begin{proof}
For any loop $u$ in $\Of$, define
$v:S^1\times[0,\ell]\to\Sigma$ by 
$$
v(s,t) = \phi^{-j}(u(s,t-j)) 
\quad\text{for}\quad (s,t)\in S^1\times[j,j+1],j<\ell.
$$
Since $v(s,0)=\phi^\ell(v(s,\ell))$ for all $s\in S^1$,
$v$ can be considered a loop in $\Omega_{\phi^\ell}$. 
Note that 
$$
-\int_{S^1\times[0,\ell]}v^*\omega = \ell\cdot\langle[\alpha_\omega],[u]\rangle.
$$
Now observe that $v(.,0)$ is freely homotopic to $\phi^{-1}(v(.,0))=v(.,1)$. 
By Corollary~\ref{cor:dehn} in Appendix~\ref{ap:dehn}, we therefore know that
$v(.,0)$ is freely homotopic to a loop $\gamma:S^1\to\fix(\phi^\ell)$. Hence, it
follows from the previous lemma with $\phi$ replaced by $\phi^\ell$, that 
$$
\int_{S^1\times[0,\ell]}v^*\omega = \int_{S^1}\gamma^*\omega \;=\; 0.
$$ 
Therefore $\langle[\alpha_\omega],[u]\rangle=0$, which ends the proof of the
proposition.  
\end{proof}

\section{Connecting orbits}\label{se:corbits}

The main result of this section is a separation mechanism for Floer connecting
orbits. 
Together with the topological separation of fixed points discussed in
Proposition~\ref{prop:fclass}, it allows us to compute the Floer homology of
diffeomorphisms of finite type. 
We expect, however, that the results of this section are applicable to a
larger class of surface diffeomorphisms.  

We start by introducing the setup. The notation is reminiscent of the
situation encountered in Section~\ref{se:diff}. However, it is only in the
next section where we return our attention to diffeomorphisms of finite type.

Let $\Sigma_0\subset\Sigma$ be a compact submanifold, not necessarily
connected. Let $N_0\subset\Sigma_0$ be a collar neighborhood of
$\p\Sigma_0$. On every connected component of $N_0$, we choose
coordinates $(q,p)\in[0,1]\times S^1$ such that
$\p\Sigma_0\cong\bigcup\{1\}\times S^1$. Let $\omega$ be an area form
on $\Sigma$ which is given by $\de q\wedge\de p$ on $N_0$.
Let $\Phi\in\symp(\Sigma,\omega)$, $H:\Sigma\to\R$ a smooth function and
$J=(J_t)_{t\in\R}$ such that the following holds:
\smallskip\\
(H1)
$\Sigma_0$ is $\Phi$-invariant. Moreover, $\Phi(x)=\psi_1(x)$ for all
$x\in\Sigma_0$, where $(\psi_t)_{t\in\R}$ denotes the Hamiltonian flow
generated by $H$. 
This means, that $\p_t\psi_t=X\circ\psi_t$, where the vector field $X$ is
defined by $\de H=\omega(X,\cdot)$.  
\smallskip\\
(H2)
There exists a constant $0<\delta<1/4$ such that on each connected component
of $N_0$, we have $\Phi(q,p)=(q,p\mp\delta)$. The sign may depend on the component. 
\smallskip\\
(H3)
$\fix(\Phi)\cap\Sigma_0=\crit(H)\cap\Sigma_0$.
\smallskip\\
(H4)
If $\omega'$ is a $\Phi$-invariant area form such that $\omega'=\omega$ on
$\Sigma\setminus N_0$, then $\alpha_{\omega'}$ has vanishing periods
on $\Omega_\Phi$. 
\smallskip\\
(H5)
For all $t\in\R$, $J_t$ is an $\omega$-compatible complex structure
which restricts to the standard complex structure on $N_0$ with
respect to the $(q,p)$-coordinates. Moreover, $J_{t+1}=\Phi^*J_t$. 
\smallskip\\
Assuming (H1--5) we prove that 
\begin{prop}[Connecting orbits]\label{prop:corbits}
Let $x^-,x^+\in\fix(\Phi)\cap\Sigma_0$ be in the same connected component of
$\Sigma_0$.  
If $u\in\M(x^-,x^+;J,\Phi)$, then $\im u \subset \Sigma_\delta$,
where $\Sigma_\delta$ denotes the $\delta$-neighborhood of 
$\Sigma_0\setminus N_0$  with respect to any of the metrics $\omega(.,J_t.)$. 
\end{prop}
To prove this proposition we vary the symplectic form. Fix $\e>0$ 
sufficiently small and set $N_\e:=\bigcup\,[\e,1-\e]\times S^1\subset N_0$. 
For every $R>0$, let $\lambda_R:\Sigma\to\R_{>0}$ be a
$\Phi$-invariant smooth function such that 
$$
\lambda_R\equiv\begin{cases}
                 R & \text{on $N_\e$}, \\
                 1 & \text{on $\Sigma\setminus N_0$}.
               \end{cases}
$$
Set 
$$
\omega_R:=\lambda_R^2\cdot\omega 
\quad\text{and}\quad 
g_{R,t}:=\omega_R(\cdot,J_t\cdot).
$$ 
Note that $\omega_R$ is $\Phi$-invariant and $\omega_R=\omega$ on
$\Sigma\setminus N_0$. 
By (H4), we can define an action functional
$\action_R:\Omega_\Phi\to\R$ such that 
\begin{eqnarray*}
\action_R(y')-\action_R(y) &=& 
\int_0^1\alpha_{\omega_R}(u(s,\cdot))\p_su(s,\cdot)\,\de s \\
&=& \int_0^1\int_0^1\omega_R(\p_tu(s,t),\p_su(s,t))\,\de t\de s
\end{eqnarray*}
for all $y,y'\in\Omega_\Phi$ and $u:\I\times\R\to\Sigma$ with
$u(s,t)=\Phi(u(s,1+t))$ and $u(0,\cdot)=y,u(1,\cdot)=y'$. 
The following observation is crucial for the proof of the proposition.
\begin{lemma}\label{le:ar}
Let $x^-,x^+\in\fix(\Phi)\cap\Sigma_0$ be in the same connected component of
$\Sigma_0$. Then 
$$
\action_R(x^+) - \action_R(x^-) = H(x^-)- H(x^+), 
$$
for every $R>0$. 
\end{lemma}
\begin{proof}
Choose a path $\gamma:\I\to\Sigma_0\setminus N_0$ from $x^-$ to $x^+$ and 
define $h:\I\times\R\to\Sigma$ by $$
h(s,t):=\psi_{-t}(\gamma(s)).
$$
By (H1), we have $h(s,t)=\Phi(h(s,t+1))$ and furthermore, by (H3),
$h(0,.)=x^-,h(1,.)=x^+$. 
Hence,
\begin{eqnarray*}
\action_R(x^+)-\action_R(x^-) 
&=& \int_0^1\!\int_0^1\!\omega_R\big(\p_t h(s,t),\p_s h(s,t)\big)\,\de s\de t \\
&=& -\int_0^1\!\Big(\int_0^1\!\de H\big(h(s,t)\big)\p_sh(s,t)\,\de s\Big)\de t \\
&=& \int_0^1\!\big(H(h(0,t))-H(h(1,t))\big)\de t \\
&=& H(x^-) - H(x^+).
\end{eqnarray*}
In the second line we use that 
$$
\omega_R(\,.\,,\p_t h(s,t))=\omega(\,.\,,\p_t h(s,t))=\de H(h(s,t)),
$$  
since $\im(h)\subset\Sigma_0\setminus N_0$.
\end{proof}
\begin{proof}[Proof of Proposition \ref{prop:corbits}]
Let $u\in\M(x^-,x^+;J,\Phi)$, i.e. $u:\R^2\to\Sigma$ is a smooth
function satisfying 
\begin{equation}\label{eq:floer}
\left\{\begin{array}{l}
u(s,t) = \phi(u(s,t+1)), \\
\p_s u + J_t(u)\p_t u = 0, \\
\lim_{s\to\pm\infty}u(s,t) = x^\pm
\end{array}\right.
\end{equation}
and
$$
\int_{\R}\int_0^1 \omega\big(\p_tu(s,t),J_t\p_tu(s,t)\big)\,\de t\de s
<\infty.
$$
We prove that 
$$
\im(u)\cap N'_\e=\emptyset,
$$ 
where 
$N'_{\e}=\bigcup\,(\e+\delta,1/2]\times S^1\subset N_0$.
Assume by contradiction that there exist $s_1<s_2$ and $t'$ such that
$u(s,t')\in N'_\e$ for all $s\in[s_1,s_2]$. 
\smallskip\\
Denote by $g$ the standard Euclidean metric on $N_0$. Note, that if $x\in N'_\e$, then
$B_g(x,\delta)\subset N_\e$. Here $B_g(x,\delta)$ denotes the $g$-disk of
radius $\delta$ around $x$. Moreover, if $x=(q,p)$, then 
$\Phi^{-1}(x)=(q,p\pm\delta)\in\p B_g(x,\delta)$.
\smallskip\\
Hence, it follows from the first equation in \eqref{eq:floer}, that
$u(s,t'+1)\in\p B_g(u(s,t'))$ for all  
$s\in[s_1,s_2]$. Fix $s\in[s_1,s_2]$ and let $r\in(0,1]$ be such that 
$$
u(\{s\}\times[t',t'+r])\subset B_g(u(s,t'),\delta),\quad
u(s,t'+r)\in\p B_g(u(s,t'),\delta).
$$
This implies that 
$$
\delta \leq \int_{t'}^{t'+r} \,|\p_tu(s,t)|_g \,\de t.
$$
By our choice of $J_t$, we know that $g_{R,t}=R^2g$ on $N_\e$, and hence that
$R\cdot|\p_tu(s,t)|_g=|\p_tu(s,t)|_{g_{R,t}}$ for all $t\in[t',t'+r]$.
Therefore,
\begin{equation*}
\begin{split}
R\cdot\delta 
&\leq \int_{t'}^{t'+r} \,|\p_tu(s,t)|_{g_{R,t}} \,\de t \\
&\leq \int_{t'}^{t'+1} \,|\p_tu(s,t)|_{g_{R,t}} \,\de t 
\,=\, \int_0^1 \,|\p_tu(s,t)|_{g_{R,t}}\,\de t.
\end{split}
\end{equation*}
In the last step we use that $|\p_tu(s,t)|_{g_{R,t}}$ is a 1-periodic function in
$t$. By H\"{o}lder's inequality we get that
$$ 
R^2\cdot\delta^2 \leq \int_0^1 \,|\p_tu(s,t)|_{g_{R,t}}^2 \,\de t,
$$
for all $R>0$ and $s\in[s_1,s_2]$. 
Now integrate over $[s_1,s_2]$:
\begin{equation}\label{eq:ineq}
\begin{split}
R^2\cdot\delta^2\cdot|s_1-s_2| 
&\leq \int_{s_1}^{s_2}\!\int_0^1\!|\p_tu(s,t)|_{g_{R,t}}^2\,\de s\de t \\
&\leq \int_\R\int_0^1\!|\p_tu(s,t)|_{g_{R,t}}^2\,\de s\de t.
\end{split}
\end{equation}
Finally, we use the energy identity
$$
|\p_tu(s,t)|_{g_{R,t}}^2=\omega_R(\p_su(s,t),\p_tu(s,t)),
$$
which follows from the second equation in \eqref{eq:floer}.
Note, that since $u$ has finite energy with respect to $\omega_R$,
the energy identity implies that 
$$
\int_\R\int_0^1\!|\p_tu(s,t)|_{g_{R,t}}^2\,\de s\de t
= 
\action_R(x^-)-\action_R(x^+).
$$
From \eqref{eq:ineq} and Lemma~\ref{le:ar} we therefore get
$$
R^2\cdot\delta^2\cdot|s_1-s_2| 
\leq
H(x^+)-H(x^-),
$$
for all $R>0$. For large $R$ this is a
contradiction and proves that $\im(u)$ is disjoint from $N'_\e$. 
Since $\e$ can be chosen arbitrarily small by an appropriate choice of
functions $\lambda_R$, $\im(u)$ is disjoint from $\bigcup\,(\delta,1/2]\times
S^1$. This proves the proposition.  
\end{proof}
\begin{remark} 
The advantage of our approach to the connecting orbits proposition compared
to Seidel's original approach in 
\cite{S1} is that we do not have to make the function $H$ depend on $R$.
Moreover, the bubbling argument completely disappears. This is relevant in the
case where $\Phi$ twists with different signs at different ``ends'' of
$\Sigma_0$. In that case, when the ends are stretched as in \cite[Lemma 4]{S1},
the energy difference of certain fixed points may go to infinity and the
bubbling argument fails.
\end{remark}

\section{Floer homology of finite type diffeomorphisms}\label{se:main}

In this section we prove Theorem~\ref{thm:main1}. 
We return to the notation of Section~\ref{se:diff}, i.e.
$\phi$ is a diffeomorphism of finite type and $\Sigma_0$ denotes the union of
connected components of $\Sigma\setminus\text{int}(N)$, where $\phi$ restricts
to the identity. 
\smallskip\\
(Monotonicity) By Propositions~\ref{prop:monotone3} and
\ref{prop:monotone4} we can choose an area form $\omega$ such that
$\phi\in\symp^m(\Sigma,\omega)$. 
Note that every Hamiltonian perturbation of $\phi$ is also in
$\symp^m(\Sigma,\omega)$.  
We impose an additional condition on $\omega$ in the next paragraph. 
\smallskip\\
(Hamiltonian perturbation) As a preparation,  
let $f_1,f_2:[0,3]\to\R$ be two functions which are constant on $[0,1]$ and
such that $f_1(q)=f_2(q)$ for all $q\in[2,3]$. Define the function 
$$
h(q) 
=\int_q^3\big(f_1(r)-f_2(r)\big)\de r,
$$
for $q\in[0,3]$. It follows that 
$$
h(q)=\begin{cases}
       \delta\cdot q + c & \text{if $q\in[0,1]$}, \\
       0                 & \text{if $q\in[2,3]$},
     \end{cases}
$$
where $\delta=f_2(0)-f_1(0)$ and $c=\int_0^3(f_1(r)-f_2(r))\de r$. 
Now consider $h(q)$ as a function on $[0,3]\times S^1$ with coordinates $(q,p)$.
The Hamiltonian vector field of $h$ with respect to $\de q\wedge\de p$ is
simply $(f_1(q)-f_2(q))\cdot\p/\p p$ at the point $(q,p)$. The time-1-maps
of the flow is thus the twist map $(q,p)\mapsto(q,p+f_1(q)-f_2(q))$.  
In particular, $(q,p)\mapsto(q,p-\delta)$ if $q\in[0,1]$.

This has the following application. Let $\Nn\subset\Sigma$ be a
$\phi$-invariant closed tubular neighborhood of $\p\Sigma_0$. On every
connected component of $\Nn$, we choose coordinates
$(q,p)\in[0,3]\times S^1$ such that $\phi$ is given by
$(q,p)\mapsto(q,p\mp f(q))$, for some monotone increasing
$f:[0,3]\to[0,1)$. Moreover, we assume that
$N_0:=\Sigma_0\cap\Nn\cong\bigcup[0,1]\times S^1$. 
Note that $f|[0,1]\equiv0$. Furthermore note, that we can 
assume that $\omega=\de q\wedge\de p$ on $\Nn$. It now follows from the
preliminary remarks, that for every $0<\delta<1/4$ there exists
$h:\Nn\to\R$ such that 
\smallskip\\
(i) 
On a connected component of $N_0$, $h(q,p)=\pm\delta\cdot q+c$. Here, the sign and
the constant $c$ may depend on the component. 
\smallskip\\
(ii)
$h\equiv0$ on $\bigcup[2,3]\times S^1$.
\smallskip\\
(iii)
Let $\psi$ denote the time-1-map of the Hamiltonian flow generated by $h$ with
respect to $\de q\wedge\de p$. For every connected component of $N_0$, there
exists a monotone increasing function $g:[0,3]\to[\delta,1)$, such that
$\phi\circ\psi(q,p)=(q,p\mp g(q))$. 
\smallskip\\
As a consequence of (i) and (ii), there exists a function $H:\Sigma\to\R$ with
$$
H(x)= \begin{cases}
     h(x) & \text{if $x\in\Nn$}, \\
     0 & \text{if $x\in\Sigma\setminus(\Sigma_0\cup\Nn)$},
    \end{cases}
$$
and such that $H|\text{int}(\Sigma_0)$ is a Morse function, 
meaning that all the critical points are non-degenerate.
We refer to \cite[Lemma 4.15]{Sc} for the extension of Morse functions.

Let $(\psi_t)_{t\in\R}$ denote the Hamiltonian flow generated by $H$ with
respect to the fixed area form $\omega$ and set
$$
\Phi:=\phi\circ\psi_1.
$$  
By construction, $\omega,\Phi,H$ and $N_0$ satisfy the hypothesis (H1,2) of the
last section. Moreover, by choosing $H$ and $\delta$ such that the $C^2$-norm
of $H$ is sufficiently small, we also guarantee (H3). 
\smallskip\\
(Fixed points) 
By (iii), $\Phi|\Nn$ has no fixed points. Since 
$\Phi=\phi$ on $\Sigma\setminus(\Sigma_0\cup\Nn)$, we therefore have
$$
\fix(\Phi) = \big(\crit(H)\cap\Sigma_0\big)\cup\big(\fix(\phi)\setminus\Sigma_0\big).
$$
In particular, $\Phi$ only has non-degenerate fixed points
and the $\Z_2$-degree of a fixed point is given by 
$$
\deg(y)=\mind_H(y)\bmod2
\quad\forall\; y\in\crit(H)\cap\Sigma_0,
$$ 
and 
$$
\deg(y)=0\bmod2
\quad\forall\; y\in\fix(\phi)\setminus\Sigma_0.
$$
The first equality follows from \cite[Lemma 7.2]{SZ}, the second from 
Proposition~\ref{prop:fclass}.
Moreover, Proposition~\ref{prop:fclass} implies that every
$y\in\fix(\phi)\setminus\Sigma_0$ forms a different fixed point class of
$\Phi$. 
This has an immediate consequence for the Floer complex $(\CF(\Phi),\p_{J})$
with respect to a generic $J=(J_t)_{t\in\R}$. 
\begin{lemma}\label{prop1}
$(\CF(\Phi),\p_{J})$ splits into the subcomplexes $(\Cc_1,\p_1)$ and
$(\Cc_2,\p_2)$, where $\Cc_1$ is generated by $\crit(H)\cap\Sigma_0$ 
and $\Cc_2$ by $\fix(\phi)\setminus\Sigma_0$. Moreover, $\Cc_2$ is graded by~$0$
and $\p_2=0$. The splitting is respected by the quantum cap action
\eqref{eq:qca}.  
\end{lemma}
\begin{proof}
If $y^\pm\in\fix(\Phi)$ are in different fixed point classes, then 
$\M(y^-,y^+;J,\phi)=\emptyset$. This follows from the first equation in
\eqref{eq:corbit} and proves the lemma.
\end{proof}
Next we show that hypothesis (H4) holds. 
\begin{lemma}[Action 2]\label{le:action2}
If $\omega'$ is a $\Phi$-invariant area form such that
$\omega'=\omega$ on $\Sigma\setminus N_0$, then $\alpha_{\omega'}$ has
vanishing periods.  
\end{lemma}
\begin{proof}
The proof is an extension of the proof of Proposition~\ref{prop:action}.
Let $u$ be a loop in $\Omega_{\Phi}$. We claim that
$\langle[\alpha_{\omega'}],[u]\rangle=0$. 
Since $\Phi^\ell$ is isotopic to $\phi^\ell$,
we can assume that $u(.,0)$ is contained either in  
$\Sigma\setminus(\Sigma_0\cup \Nn)$ or in $\Sigma_0\setminus N_0$.
The argument is similar as in the above mentioned proof and relies on Lemma
~\ref{lemma:mu} and Corollary~\ref{cor:dehn}. 
The first case reduces to the case considered in
Proposition~\ref{prop:action}.
In the second case, define $h:S^1\times\R\to\Sigma_0\setminus N_0$ by 
$h(s,t):=\psi_{-t}(u(s,0))$; $h$ is a loop in $\Omega_\Phi$. 
Since $\omega'=\omega$ on $\im(h)$, it follows that 
$\int h^*\omega'=\int h^*\omega=0$ and
hence, from Lemma~\ref{lemma:mu}, that $\int u^*\omega'=0$. 
\end{proof}
\noindent
(Path of complex structures)
Let $J_0$ be a $\omega$-compatible complex structure on $\Sigma$ which
restricts to the standard complex structure on $N_0$. 
Let $J=(J_t)_{t\in\R}$ be a smooth path of $\omega$-compatible complex
structures such that $J_{t+1}=\Phi^*J_t$ and
$J_t(x)=(\psi_t^*J_0)(x)$ for all $t\in\R$ and $x\in\Sigma_0$. 
The existence of such a $J$ relies on the contractibility of the
space of $\omega$-compatible complex structures on $\Sigma$. 
Note that $J$ satisfies (H5). 
Below, we impose an additional regularity condition on $J_0$. 

We are now in position to apply Proposition~\ref{prop:corbits} and compute the
homology of $(\Cc_1,\p_1)$. Let $\p_{\pm}\Sigma_0$ denote the union of
components of $\p\Sigma_0$ where in a neighborhood, $\Phi$ is given by
$(q,p)\mapsto(q,p\mp\delta)$. 
\begin{lemma}\label{prop2}
The homology of $(\Cc_1,\p_1)$ is isomorphic to 
$H_*(\Sigma_0,\p_+\Sigma_0;\Z_2)$. 
The quantum cap product is given by the ordinary cap product 
$$
H^*(\Sigma;\Z_2)\otimes H_*(\Sigma_0,\p_+\Sigma_0;\Z_2)
\To H_*(\Sigma_0,\p_+\Sigma_0;\Z_2).
$$ 
\end{lemma}
\begin{proof}
The proof is by the same technique as in \cite{S1}. 
By modifying $J_0$ in a neighborhood of $\crit(H)\cap\Sigma_0$, we can
assume that $\nabla H$ is a Morse-Smale vector field on $\Sigma_0$, where 
the gradient is with respect to the metric $\omega(.,J_0.)$.
This means that stable and unstable manifolds are all transverse to each other
and that the stable, unstable manifolds are all transverse to $\p\Sigma_0$. 

Note that by Proposition~\ref{prop:fclass}, $(\Cc_1,\p_1)$ splits into
subcomplexes generated by fixed points of $\Phi$ which are in the same
connected component of $\Sigma_0$. Let $x^\pm$ be a pair of such fixed points
and set $\M:=\M(x^-,x^+;J,\Phi)$.
For every $u\in\M$, it follows from Proposition~\ref{prop:corbits} that
$\im(u)\subset\Sigma_0$.  
Define the map $\tilde{u}:\R^2\to\Sigma_0,(s,t)\mapsto\psi_t(u(s,t))$. 
A straight forward calculation, using that
$u(s,t)=\psi_1(u(s,t+1)),\psi_t\circ\psi_1=\psi_{t+1}$ and 
$J_t=\psi_t^*J_0$ on $\im(u)$, shows that
\begin{equation}\label{eq:morse}
\left\{\begin{array}{l}
\tilde{u}(s,t) = \tilde{u}(s,t+1), \\
\p_s\tilde{u}+J_0(\tilde{u})\big(\p_t\tilde{u}-X_H(\tilde{u})\big) = 0,\\
\lim_{s\to\pm\infty}\tilde{u}(s,t) = x^\pm,
\end{array}\right.
\end{equation}
where $X_H$ denotes the Hamiltonian vector field of $H$. 
The system~\eqref{eq:morse} was studied in 
\cite[Theorem 7.3]{SZ}. There it is 
shown that if $H|\Sigma_0$ is replaced by $\e\cdot H|\Sigma_0$ with $\e>0$
sufficiently small, then every solution of \eqref{eq:morse} is
independent of the $t$-variable and is therefore a solution of  
$$
\de\tilde{u}/\de s = \e\cdot\nabla H(\tilde{u}).
$$
Moreover, $\tilde{u}$ is regular in the sense that the operator which is
defined by linearizing the equations~\eqref{eq:morse} is surjective, see
page~\pageref{eq:corbit}.  
If we go back to the definition of $H$, we can assume from now on that $\e=1$. 
That the ambient space $\Sigma_0$ has non-empty boundary does not affect the
argument in \cite[Theorem 7.3]{SZ}. It is essential however, that
$\pi_2(\Sigma_0)=0$.  

It follows that every $u\in\M$ is regular and that the map $u\mapsto\tilde{u}$
induces a diffeomorphism of $\M$ and the space of (parameterized) flow lines of
$\nabla H$ which are contained in $\Sigma_0$ and connect the critical points $x^\pm$. 
Furthermore, these diffeomorphisms, one for each pair $x^\pm$, induce an isomorphism
$(\Cc_1,\p_1)\cong(\CM(H|\Sigma_0),\p_{\nabla\!H})$ of chain complexes. Here,
$\CM(H|\Sigma_0)$ is freely 
generated by $\crit(H)\cap\Sigma_0$
and $\p_{\nabla\!H}$ is defined by counting index-1 flow lines of $\nabla H$
which are contained in $\Sigma_0$. 
Note that $\nabla H$ points outwards/inwards at a
component of $\p_+/\p_-\Sigma_0$. 
The homology of $(\CM(H|\Sigma_0),\p_{\nabla\!H})$ is therefore isomorphic to 
$H_*(\Sigma_0,\p_+\Sigma_0;\Z_2)$. See \cite{Sc} for details on relative Morse
homology. 

Similarly we can identify the quantum cap product, defined in \eqref{eq:qca}.   
Note that the image of the evaluation
map $\M\to\Sigma,u\mapsto u(0,0)$, is $W^u(\nabla H,x^-)\cap W^s(\nabla H,x^+)$. 
Choose a Morse function $f:\Sigma\to\R$ such that the evaluation map is
transverse to $W^u(\nabla f,x)$ for all $x\in\crit(f)$. 
For $x\in\crit(f)$ and $x^\pm\in\crit(H)\cap\Sigma_0$ with
$\mind_H(x^+)=\mind_H(x^-)+\mind_f(x)$, let 
$q(x;x^-,x^+)\in\Z_2$ be the cardinality mod 2 of
$W^u(\nabla f,x)\cap W^u(\nabla H,x^-)\cap W^s(\nabla H,x^+)$. 
The map 
$$
\CoM(f)\otimes \CM(H|\Sigma_0)\To \CM(H|\Sigma_0),\quad
x\otimes y\Mapsto\sum_{z}q(x;y,z)z,
$$
which induces the quantum cap product on homology, is therefore given in
purely Morse theoretical terms. On the level of homology, it is the ordinary
cap product. This finishes the proof of the lemma.
\end{proof}
\begin{proof}[Proof of Theorem \ref{thm:main1}]
From Lemma~\ref{prop1} and \ref{prop2}, it follows that  
$$
\HF(\phi) \cong
H_*(\Sigma_0,\p_+\Sigma_0;\Z_2)\oplus
\Z_2^{\#\fix(\phi|\Sigma\setminus\Sigma_0)}.
$$
Moreover, $H^*(\Sigma;\Z_2)$ acts on the first summand by ordinary cap
product.
Since every fixed point of $\phi|\Sigma\setminus\Sigma_0$ has fixed point
index~1, the Lefschetz fixed point formula implies that 
$$
\#(\fix(\phi)\setminus\Sigma_0) = \Lambda(\phi|\Sigma\setminus\Sigma_0).
$$
It remains to show that $1\in H^0(\Sigma;\Z_2)$ acts on
$\Z_2^{\Lambda(\phi|\Sigma\setminus\Sigma_0)}$ by
the identity and any element of $H^1(\Sigma;\Z_2)\oplus H^2(\Sigma;\Z_2)$ by
the zero map. 

From the proof of Lemma~\ref{prop1}, we know that if
$\M(y^-,y^+;J,\Phi)\ne\emptyset$ for some
$y^-,y^+\in\fix(\phi)\setminus\Sigma_0$, then $y^-=y^+$.
Since the action $\alpha_\omega$ has vanishing periods on $\Omega_\Phi$, by
Lemma~\ref{le:action2}, it follows that every $u\in\M(y,y;J,\Phi)$ has energy
zero.
Hence, $\M(y,y;J,\Phi)=\M_0(y,y;J,\Phi)$ consists of the constant map. 
Now choose a Morse function $f:\Sigma\to\R$ with only one critical point $x_0$
of index 0 and such that 
$\fix(\phi|\Sigma\setminus\Sigma_0)\subset W^u(\nabla f,x_0)$. 
It follows that $q(x;y^-,y^+)\ne0$ if and only if $x=x_0$ and $y^-=y^+$. 
This ends the proof.
\end{proof}
\begin{remark}
Using Theorem~\ref{thm:main1}, we can verify Seidel's result~\cite[Theorem
1]{S2} in the special case of an algebraically finite mapping class: if $g$ is
a non-trivial mapping class,  then the quantum cap product is trivial on
$H^2(\Sigma;\Z_2)\otimes\HF(g)$. 

To see this, assume that $g$ is algebraically finite.
By Theorem~\ref{thm:main1}, the only possibly non-trivial part of the the
quantum cap product is given by the cap product on
$H^*(\Sigma;\Z_2)\otimes H_*(\Sigma_0,\p_+\Sigma_0;\Z_2)$.
The submanifold $\Sigma_0\subset\Sigma$ has non-trivial boundary since 
$g$ is non-trivial. 
Now the cap product on 
$H^*(\Sigma;\Z_2)\otimes H_*(\Sigma_0,\p_+\Sigma_0;\Z_2)$
factors through the homomorphism 
$\iota^*:H^*(\Sigma;\Z_2)\to H^*(\Sigma_0;\Z_2)$ which is induced by the
inclusion $\iota:\Sigma_0\inclusion\Sigma$. 
Since $H^2(\Sigma_0;\Z_2)=0$, this proves the claim.
 
Similarly, it follows from Theorem~\ref{thm:main1}, that if 
$\alpha\in H^1(\Sigma;\Z_2)$ acts 
non-trivially on $\HF(g)$, then there exists a map $\gamma:S^1\to\Sigma_0$
such that $\langle\alpha,[\gamma]\rangle=1$. This is a special case of
\cite[Theorem 2]{S2}.  
\end{remark}

\section{Isolated plane curve singularities}\label{se:singularity}

In this section we prove Theorem~\ref{thm:main2} and Corollary~\ref{cor:dehn1}.
We begin with a brief summary of the basic facts on isolated plane curve
singularities. 
The standard reference is Milnor's book \cite{Mi1}. 

An {\bf isolated plane curve singularity} is a germ $[f]$ of holomorphic functions
$f:(U,0)\to(\C,0)$, where $U\subset\C^2$ is a neighborhood of $0$, 
with $(\de f)^{-1}(0)=\{0\}$. Let $f:U\to\C$ be such a function. 
For $\e>0$ sufficiently small, the singular fiber $f^{-1}(0)$ intersects
transversally with the 3-sphere $S_\e:=\{(x,y)\in\C^2:|x|^2+|y|^2=\e\}$. 
The intersection $L:=f^{-1}(0)\cap S_\e\subset S_\e$ is a compact oriented
1-manifold, i.e. a link. A link obtained in this way is called an 
{\bf algebraic link}. 
An algebraic link is a fibred link: the map 
$$
\pi:S_\e\setminus L\To \{z\in\C:|z|=1\},
\quad z\Mapsto f(z)/|f(z)|,
$$
is a fibration, the Milnor fibration, and the {\bf Milnor fiber}
$M:=\pi^{-1}(1)\cup L$ 
is a Seifert surface of $L$. This means that $M\subset S_\e$ is an 
compact connected oriented embedded 2-manifold with $\p M=L$. 

The {\bf geometric monodromy} is an isotopy class of the group
$\diff^+_c(M)$ of orientation preserving diffeomorphisms which are the
identity near $\p M$ and is defined as follows.
Given a connection on $S_\e\setminus L$, i.e. a rank-1 subbundle of the
tangent bundle of $S_\e\setminus L$ which is transversal to $\ker(\de\pi)$, 
parallel transport induces an orientation preserving diffeomorphism of
$\pi^{-1}(1)$; a so-called characteristic diffeomorphism. 
To extend this diffeomorphism to $M$, we specify the connection in a
neighborhood of $L$. For this observe the following.

Let $L'$ be a connected component of $L$ and $T\subset S_\e$ be a
tubular neighborhood of $L'$. 
A standard meridian of $(T,L')$ is an embedded circle in 
$T\setminus L'$, which is homologically trivial in $T$ and has linking 
number $1$ with $L'$ in $S_\e$. There is a fibration of $T\setminus L'$ such that
every fiber is a standard meridian of $(T,L')$.
This fibration is unique up to isotopy and if $T$ is sufficiently
small, induces a connection on $T\setminus L'$. 
In this way, we get a standard connection in neighborhood of $L$. 
A connection on $S_\e\setminus L$ which restricts to this
standard connection induces a characteristic
diffeomorphism which is compactly supported and hence extends trivially to $M$.
Moreover, any two such diffeomorphisms are isotopic in $\diff^+_c(M)$.
The isotopy class obtained in this way is therefore an invariant of $\pi$. 

The proof of Theorem~\ref{thm:main2} relies on the following result, which is
a refinement of the classical result of A'Campo~\cite{AC2} and L{\^e}~\cite{L}
that the Lefschetz number of the geometric monodromy of an isolated plane
curve singularity vanishes.   
We use the following notation: 
we denote by $\iota:\diff^+_c(M)\to\diff^+(M,\p M)$ the inclusion, where 
$\diff^+(M,\p M)$ denotes the group of orientation preserving
diffeomorphisms which are the identity on $\p M$. 
\begin{prop}\label{prop:monodromy}
Let $M$ be the Milnor fiber and $g$ be the geometric monodromy of an isolated
plane curve singularity.  
There is a representative $\phi\in\iota_*g$ which is a diffeomorphism of
finite type (same definition as for closed surfaces) and such that
$\fix(\phi)=\p M$.
Moreover, $\phi$ only has positive twists. 
\end{prop}
As already mentioned in the introduction, this proposition follows from
the work of A'Campo~\cite{AC1}, \cite{AC4} on the geometric monodromy. 
Our proof, given in
Appendix~\ref{ap:monodromy}, relies on the work of Eisenbud and
Neumann~\cite{EN} on the monodromy of plane curve singularities.  
The relevant results from \cite{EN} are summarized in Appendix~\ref{ap:monodromy}.
We remark, that the vanishing of the Lefschetz number of diffeomorphism of
finite type is not enough to exclude
the existence of pointwise fixed annuli. 
\begin{proof}[Proof of Theorem~\ref{thm:main2}]
We recall that $\Sigma$ denotes a closed oriented surface of genus $\geq2$ and 
$M\subset\Sigma$ the Milnor fiber of an isolated plane curve
singularity.
Moreover, $g$ denotes the mapping class of $\Sigma$ which is obtained by
extending the geometric monodromy of the singularity trivially to $\Sigma$.

Assume for the moment that no component of $\Sigma\setminus\text{int}(M)$ is a disk.
From Proposition~\ref{prop:monodromy}, it follows that there exists a
representative $\phi\in g$ which is of finite type and such that 
$\Sigma_0=\Sigma\setminus\text{int}(M)$, $\p_+\Sigma_0=\p M$ and 
$\Lambda(\phi|\text{int}(M))=0$. 
Theorem~\ref{thm:main1} therefore implies that 
$$
\HF(\phi)\cong H_*(\Sigma\setminus\text{int}(M),\p M;\Z_2).
$$
By excision, it follows that $\HF(\phi)\cong H_*(\Sigma,M;\Z_2)$.

Now assume that $D_1,\dots,D_n$ are disk components of
$\Sigma\setminus\text{int}(M)$. 
Set $M':=M\cup D_1\cup\cdots\cup D_n$.
We claim that there exists a representative $\phi\in g$ which is of finite
type and such that 
$\Sigma_0=\Sigma\setminus\text{int}(M')$, $\p_+\Sigma_0=\p M'$ and 
$\Lambda(\phi|\mathrm{int}(M'))=n$. Theorem~\ref{thm:main1} then implies
that 
$$
\HF(\phi)
\cong H_*(\Sigma\setminus\text{int}(M'),\p M';\Z_2)\oplus\Z_2^n
\cong H_*(\Sigma,M';\Z_2)\oplus\Z_2^n.
$$
Since $H_*(\Sigma,M';\Z_2)\oplus\Z_2^n\cong H_*(\Sigma,M;\Z_2)$,
this proves Theorem~\ref{thm:main2} up to the claim above.
The claim follows from Proposition~\ref{prop:monodromy} by ``collapsing'' 
each disk component of $\Sigma\setminus\text{int}(M)$ to a point.  
\end{proof}
Next we prove Corollary~\ref{cor:dehn1}. First, we recall some terminology.
Let $k\in\N_{>0}$. An $A_k${\bf -configuration} in $\Sigma$ is a $k$-tuple
$(C_1,\dots,C_k)$ of embedded circles in $\Sigma$ such that 
$$
\#(C_{i}\pitchfork C_{i+1})=1 
\quad\text{for } i=1,\dots,k-1, 
\quad
\#(C_i\pitchfork C_{j})=0
\quad\text{if } |i-j|>1.
$$ 
The $A_k${\bf -singularity} is the germ of the function $f(x,y)=x^2+y^{k+1}$. 
It is a classical result in the theory of singularities,
that the Milnor fiber $M$ of this singularity
contains an $A_k$-configuration, which is (i) a spine of $M$ and (ii) a
distinguished basis of vanishing cycles. See \cite[Section 2.9]{Ar}, \cite{AC3}.
\begin{proof}[Proof of Corollary \ref{cor:dehn1}] 
Let $(C_1,\dots,C_k)$ be an $A_k$-configuration in $\Sigma$. 
From the remark (i) above, it follows that a tubular neighborhood $N$ of
$C_1\cup\cdots\cup C_k$ can be identified with the Milnor fiber of the
$A_k$-singularity such that, by (ii), $C_1\cup\cdots\cup C_k$ is a distinguished
set of vanishing cycles. This implies, that the class $g$ of the
product $\tau_1\circ\cdots\circ\tau_k$ of right Dehn twists can also be obtained
by extending the geometric 
monodromy of the $A_k$-singularity trivially to $\Sigma$. It
therefore follows from Theorem~\ref{thm:main2} that 
$$
\HF(g) \cong H_*(\Sigma,N;\Z_2) 
\cong H_*(\Sigma,C_1\cup\cdots\cup C_k;\Z_2).
$$
This together with naturality of Floer homology proves the corollary.
\end{proof}
\begin{proof}[Proof of Theorem \ref{thm:main3}]
Let $M$ be the Milnor fiber and $g$ be the geometric monodromy of an isolated
plane curve singularity. 
Let $\phi\in\iota_*g$ be as in Proposition~\ref{prop:monodromy}. 
Since $\phi$ only has positive twists, we can perturb it near the
boundary to a diffeomorphism $\phi_+$ such that $\fix(\phi_+)=\emptyset$. 
Furthermore, it follows as in the proof of Proposition~\ref{prop:monotone3},
that if $\omega$ is a $\phi$-invariant area form on $M$, then
$[\omega_\phi]=0$. Hence, $\phi$ is monotone in the sense defined in
Appendix~\ref{ap:open}. Therefore, $\HF(g,+)=\HF(\phi,+)=0$.   
\end{proof}

\appendix

\section{Products of disjoint Dehn twists}\label{ap:dehn}

The goal of this appendix is to prove Proposition~\ref{prop:dehn}, 
which was also stated in \cite[Lemma 3]{S1}.
This result is used for the proof of the fixed point as well as the action
proposition in section\ref{se:diff}.
Let $I$ denote either $\I$ or $S^1$. For every continuous map 
$\gamma:I\to\Sigma$ and compact 1-dimensional submanifold $C\subset\Sigma$, 
there is the geometric intersection number
$$
i(\gamma,C) = 
\min\{\#\beta^{-1}(C)\,|\,\beta \text{ is homotopic rel endpoints to } \gamma\}.
$$
If $I=S^1$, homotopic rel endpoints means freely homotopic. 
\begin{prop}\label{prop:dehn}
Let $C\subset\Sigma$ be a finite union of non-contractible disjoint circles. 
Let $\phi$ be a product of Dehn twists along $C$ which twists with the same
sign along parallel components of $C$.  
Let $\gamma:I\to\Sigma$ be such that $\gamma(\p I)\cap\supp(\phi)=\emptyset$.
If $\gamma$ is homotopic rel endpoints to $\phi\circ\gamma$, 
then $i(\gamma,C)=0$. 
\end{prop}
An immediate consequence is the following
\begin{cor}\label{cor:dehn}
Let $\phi$ be of finite type and $\ell>0$ be as such that
$\phi^\ell|\Sigma\setminus N=\id$. 
Let $\gamma:I\to\Sigma$ be such that $\gamma(\p I)\cap\supp(\phi^\ell)=\emptyset$.
If $\gamma$ is homotopic rel endpoints
to $\phi^\ell\circ\gamma$, then there exists $\gamma':I\to\fix(\phi^\ell)$ which is
homotopic rel endpoints to $\gamma$.
\end{cor}
\begin{proof}[Proof of Proposition \ref{prop:dehn}]
We begin with the case $I=\I$. For every component of $C$, let $N'$ be the
closed tubular neighborhood, where the Dehn twist is supported.
For different components of $C$, these tubular neighborhoods
are disjoint. Let $N$ be the union of the $N'$ and $\Ce:=\p N$. We prove in
several steps that $i(\gamma,\Ce)=0$, if $\gamma$ is homotopic rel endpoints to
$\phi\circ\gamma$. 

Let $\beta:\I\to\Sigma$  and $u:\I^2\to\Sigma$ be such that
$\#\beta^{-1}(\Ce)=i(\gamma,\Ce)$ and 
$$ 
u(s,0)=\beta(s),\quad u(s,1)=\phi(\beta(s)),\quad 
u(0,t)=\beta(0), \quad u(1,t)=\beta(1)
$$
for all $s,t\in\I$. 
We assume that $\beta(0),\beta(1)\in\Sigma\setminus N$. 
Without loss of generality we further assume that $u$ is
transverse to $\Ce$. Hence $B:=u^{-1}(\Ce)\subset \I^2$ is a compact
1-dimensional submanifold with boundary 
$\p B=\beta^{-1}(\Ce)\cup (\phi\circ\beta)^{-1}(\Ce)$. 
Every component of $B$ is either a circle or an arc. 
\begin{claim}{1}
We may assume that no component of $B$ is a circle.
\end{claim}
\begin{pc}
Let $S$ be a circle component of $B$. The interior of $S$, 
i.e. the region of $\I^2$ bounded by $S$, is a disk. The restriction of $u$ to
this disk induces an element of 
$\pi_2(\Sigma,C')$, where $C'$ is the component of $\Ce$ where $S$ is mapped to.
Since $C'$ is non-contractible, $\pi_2(\Sigma,C')=0$, and hence $u$ can be
deformed in a neighborhood of the interior of $S$ in such a way that $B$ has
less circle components. Repeating this argument finitely many times proves
claim 1.
From now on, we assume that every component of $B$ is an arc.
\end{pc}
\begin{claim}{2}
There is no component $B'$ of $B$ with $\p B'\subset\I\times0$ or 
$\p B'\subset\I\times1$.
\end{claim}
\begin{pc}
Assume that $B'$ is such that $\p B'\subset\I\times0$. Note that the
boundary points of $B'$ are intersection points of $\beta$ with $\Ce$. Since $B'$
is mapped to some component of $\Ce$ under $u$, it follows that these
intersection points can be removed by a homotopy of $\beta$. This however, contradicts
the definition of $\beta$. Hence every $B'$ with one boundary point on $\I\times0$
has the other boundary point on $\I\times1$. On the other hand, $B$ has the
same the number of boundary points on $\I\times0$ as on $\I\times1$, since
$\beta$ and $\phi\circ\beta$ intersect $\Ce$ in the same number of points. This
proves claim 2. The rest of the proof is devoted to 
\end{pc}
\begin{claim}{3}
$B$ is empty, i.e. $\#(\im\beta\cap\Ce)=0$.
\end{claim}
\begin{pc}
The proof is by contradiction; assume that $B\ne\emptyset$.
First we define an integer valued function on $\pi_0(B)$. Let $B'$ be a
component of $B$. By claim 1 and 2, $B'$ is an arc and $\p B'=\{(b,0),(b,1)\}$,
for some $b\in(0,1)$. Let $C'$ be the component of $\Ce$ where $B'$ is mapped
to under $u$. Since $u(b,1)=\phi(u(b,0))=u(b,0)$, the restriction $u|B'$ has a
well defined mapping degree $\deg u|B'$, once orientations of $B'$ and $C'$
are fixed. We orient $B'$ ``from the bottom to the top'' and
choose an orientation of $\Ce$ such that the orientations of two homotopic
components match. We thus have the map  
$$
\de:\pi_0(B)\To\Z,\quad B'\Mapsto\deg u|B'.
$$
Let $\pi_0(B)=\{B_1,\dots,B_{2n}\}$ be ordered such that 
$$
b_i<b_j \quad\Iff\quad i<j,
$$
where $\p B_i=\{(b_i,0),(b_i,1)\}$. Note that the cardinality of $\pi_0(B)$ is
even, since $\beta(0)$ and $\beta(1)$ are both in the complement of $N$. 
Observe that the loops $u|B_1$ and $u|B_{2n}$ are both homotopic to the
constant loop and hence 
$$
\de(B_1) = \de(B_{2n}) = 0.
$$
To prove claim 3, we will now show that $\de(B_{2n})\ne0$, which is a
contradiction. 
In fact, we will prove by induction that 
\begin{equation}\label{eq:ind}\tag{$\star_k$}
\de(B_{2k})\ne0, \quad
\sign(b_{2k}) = \sign(b_{2k-2}) \quad\mathrm{and}\quad
C_{2k}\sim C_{2k-2},
\end{equation}
for all $1\leq k\leq n$.
Here $\sign(b_i)$ denotes the sign of $b_i$ as an intersection point of
$\beta$ and $\Ce$ and $C_i$ is the component of $\Ce$ where $B_i$ is mapped to
under $u$. We will use the following formula for the function $\de$, which
we prove later on:
\begin{equation}\label{eq:de}
\de(B_{2k}) = -\sum_{i=1}^k\sign(b_{2i})\cdot\e_{2i},
\end{equation} 
for all $1\leq k\leq\ell$. Here $\e_i=\pm$ is the sign of the twist of $\phi$ 
along $C_i$. For $k=1$, formula (\ref{eq:de}) gives
$\de(B_2)=\sign(b_2)\cdot\e_2\ne0$. This proves the induction
hypothesis $(\star_1)$, since the other two statements are empty in this case.

Assume that $(\star_k)$ holds for all $1\leq i\leq k<n$. To show that it
also holds for $k+1$ first note, that if 
\begin{equation}\label{eq:ind1}
\sign(b_{2k})=\sign(b_{2k+2}) \quad\mathrm{and}\quad C_{2k}\sim
C_{2k+2}, 
\end{equation}
then by formula \eqref{eq:de}, we get that
$$
\de(B_{2k+2}) \,=\, 
(k+1)\cdot\sign(b_{2k+2})\cdot\e_{2k+2} \,\ne\, 0.
$$ 
To prove \eqref{eq:ind1}, we consider the restriction of $u$ to the region
$D\subset \I^2$ that is bounded by the arcs 
$B_{2k},B_{2k+1}$ and $[b_{2k},b_{2k+1}]\times0,[b_{2k},b_{2k+1}]\times1$.
Note that since $u(D)\subset\Sigma\setminus\text{int}(N)$, $u(s,1)=u(s,0)$ for all
$s\in[b_{2k},b_{2k+1}]$. Hence $u|D$ gives a homotopy between the loops
$\de(B_{2k})\cdot C_{2k}$ and $\de(B_{2k+1})\cdot C_{2k+1}$. This implies that
$C_{2k}\sim C_{2k+1}$ and that $\de(B_{2k})=\de(B_{2k+1})$. 
Since $C_{2k+1}$ and $C_{2k+2}$ are both contained in the boundary of one
component of $N$, they are clearly homotopic and hence $C_{2k}\sim C_{2k+2}$. 
It thus remains to show that $\sign(b_{2k})=\sign(b_{2k+2})$. For this we
need 
\begin{claim}{4}
$C_{2k}\ne C_{2k+1}$ and $C_{2k+1}\ne C_{2k+2}$.
\end{claim}
\begin{pc}
Assume that $C_{2k}=C_{2k+1}$ and consider the path 
$s\mapsto u(s,0),s\in[b_{2k},b_{2k+1}]$. 
We claim that since $\de(B_{2k})=\de(B_{2k+1})\ne0$, the path can be
deformed into $C_{2k}$, which contradicts minimality of $\beta$. 
For the proof of the claim, we refer to \cite[Corollary 5]{Ga1}.
Similarly,   
$C_{2k+1}=C_{2k+2}$ is also a contradiction to minimality of $\beta$.
\end{pc}
\begin{claim}{5}
$\sign(b_{2k})=\sign(b_{2k+1})=\sign(b_{2k+2})$.
\end{claim}
\begin{pc}
Since $C_{2k}$ and $C_{2k+1}$ are homotopic and disjoint, there exists
an embedded annulus $A\subset\Sigma$ which bounds $C_{2k}$ and
$C_{2k+1}$. Since
$\de(B_{2k})=\de(B_{2k+1})\ne0$, it follows that $u(D)=A$, where $D$ is as
above. 
Otherwise, we could find a map from the torus to $\Sigma$ with nonzero
degree. Such a map does not exist since the genus of $\Sigma$ is $>1$. Thus we
know that the path $\beta$ enters $A$ at $b_{2k}$ and leaves
at $b_{2k+1}$. The intersection points therefore have the same sign,
$\sign(b_{2k})=\sign(b_{2k+1})$. 

The proof that $\sign(b_{2k+1})=\sign(b_{2k+2})$ is similar. In this case
however, we know that the region $D'$ of $\I^2$ between $B_{2k+1}$ and
$B_{2k+2}$ is mapped to a component $N'$ of $N$ under $u$. Since
$\de(B_{2k+1})\ne0$, it follows that $u(D')=N'$ and hence that $b_{2k+1}$ and
$b_{2k+2}$ have the same sign, as above. This proves claim 5 and ends the
proof of \eqref{eq:ind1}. It remains to proof formula~\eqref{eq:de}. 
\end{pc}
\smallskip\\
Formula \eqref{eq:de} is the consequence the following properties of the
function $\de$. First of all, as was shown above, we have
that 
\begin{equation}\label{eq:de1}
\de(B_{2k})=\de(B_{2k+1}),
\end{equation}
for all $1\leq k\leq n-1$. We now show that
\begin{equation}\label{eq:de2}
\de(B_{2k})=\de(B_{2k-1})-\sign(b_{2k})\cdot\e_{2k},
\end{equation}
for all $2\leq k\leq\ell$. The idea is to look again at the region
$D'\subset\I^2$ which is bounded by the arcs $B_{2k-1},B_{2k}$ and
$[b_{2k-1},b_{2k}]\times0,[b_{2k-1},b_{2k}]\times1$. Let $N'$ be the component
of $N$ where $D'$ is mapped to. As in claim 4, we conclude that
$\p N'=C_{2k-1}\cup C_{2k}$. Choose an orientation preserving diffeomorphism
$[0,1]\times S^1\cong N'$ such that $C_{2k-1}\cong 0\times S^1$ and 
$C_{2k}\cong 1\times S^1$. Assume for the moment that the orientations of
$C_{2k-1}$ and $C_{2k}$ are compatible with these diffeomorphisms. This is
equivalent to saying that $\sign(b_{2k-1})=1$. Let
$\pr:[0,1]\times S^1\to S^1$ denote the projection onto the second factor. The map
$\pr\circ u:\p D\to S^1$ is the composition of four loops, denoted by $v$,
$u|B_{2k}$, the inverse of $w$ and the inverse of
$u|B_{2k-1}$. Here $v,w$ are the loops
$s\mapsto\pr(u(s,0)),s\mapsto\pr(u(s,1))$ for $s\in[b_{2k-1},b_{2k}]$. 
Since the mapping degree of $\pr\circ u|\p D$ vanishes, it follows that 
$$
\deg v + \de(B_{2k}) - \deg w - \de(B_{2k-1}) \,=\, 0.
$$ 
However, $\deg v-\deg w=\e_{2k}$, and this proves equation~\eqref{eq:de2} 
in the case that $\sign(b_{2k})=1$. 
If $\sign(b_{2k})=-1$, $\pr$ induces
orientation reversing maps on $C_{2k-1}$ and $C_{2k}$. Hence,
the same argument as above gives  
$$
\deg v - \de(B_{2k}) - \deg w + \de(B_{2k-1}) \,=\, 0,
$$ 
which again implies \eqref{eq:de2}. Formula \eqref{eq:de} is now an immediate
consequence of \eqref{eq:de1},\eqref{eq:de2} and the fact that $\de(B_1)=0$. 
This ends the proof of the proposition in the case $I=\I$. 
\end{pc}

The proof in the case $I=S^1$ follows the same line of arguments as above. 
In this case, however, $B=u^{-1}(\Ce)$ is a subset of $S^1\times\I$ instead if $\I^2$. 
This does not affect the proofs of claims 1 and 2 above and hence, we can assume that
every component $B'$ of $B$ is an arc with $\p B'=\{(b,0),(b',1)\}$. Note that
$b,b'\in S^1$ are not necessarily equal. We show how the proof of claim 3
extends to this case. The idea is to consider the $\ell$-fold catenation  
$v:S^1\times[0,\ell]\to\Sigma$, defined by 
$$
v(s,t):=\phi^j(u(s,t-j)) 
\quad\text{for}\quad (s,t)\in S^1\times[j,j+1],j<\ell,
$$
where $\ell>0$ is chosen such that
every component $A'$ of $A:=v^{-1}(\Ce)$ satisfies $\p A'=\{(a,0),(a,\ell)\}$ for
some $a\in S^1$. The rest of the argument is analogous to the above. Define
the function $\de:\pi_0(A)\to\Z$ and prove that it satisfies
\begin{equation}\label{eq:de3}
\de(A_{2k}) = \de(A_1) - k\cdot\ell\cdot\sign(a_{2k})\cdot\e_{2k}.
\end{equation}
Here, the signs $\sign(a_i)$ and $\e_i$ are defined as above. The
additional factor $\ell$ appears since $\phi^\ell$ twists with multiplicity
$\ell$ along every component of $C$. The order on $\pi_0(A)$ is a cyclic order
induced from the boundary points as above. Formula~\eqref{eq:de3} is therefore
a contradiction, if $\pi_0(A)\ne\emptyset$.  
This ends the proof of the proposition. 
\end{proof}

\section{Decomposition of the monodromy}\label{ap:monodromy}

In this appendix we prove Proposition \ref{prop:monodromy}.
For this we use the theory of splice diagrams which was developed by
Eisenbud and Neumann \cite{EN}, \cite{Ne} to compute invariants of plane curve
singularities. In the following we summarize their results on the
geometric monodromy. We would like to emphasis here, that our discussion is
far from being self-contained. Instead, we only concentrate on the facts which
are useful for our purpose. For details, proofs as well as the general picture,
we refer to the excellent monograph \cite{EN}.
Let us introduce some terminology. 
\begin{definition}
We denote by $M$ a compact oriented 2-manifold with boundary. 
\smallskip\\
(i) 
Let $\phi:M\to M$ be a diffeomorphism and $C\subset M$ be a $\phi$-invariant
union of disjoint non-contractible circles. A $\phi${\bf -component} of
$(M,C)$ is a a union of connected components of $M\setminus C$
which are cyclically permuted by $\phi$. 
The {\bf topological type} of a $\phi$-component $M'$ is the triple
$(\chi,d,h)$, where $\chi$ is the Euler characteristic, $d$ the number
of connected components and $h$ the number of ends of $\text{int}(M')$. 
\smallskip\\
(ii) 
An orientation preserving diffeomorphism $\phi:M\to M$ is called an
{\bf admissible twist map} if $M$ is a union of annuli and:
\smallskip\\
(1)
$\phi$ cyclically permutes the connected components of $M$.
\smallskip\\
(2)
If $n>0$ denotes the number of connected components of $M$, then
$\phi^n$ is given by the local model
$$
[0,1]\times S^1\ni(q,p)\Mapsto \big(q,p-f(q)\big),
$$
where $f:[0,1]\to\R$ is monotone. 
\smallskip\\
(3)
There exists $q>0$ such that $\phi^q|\p M=\id$.
\smallskip\\
If $\phi:M\to M$ is an admissible twist map, then its {\bf twist number}
is defined by  
$$
\ell:=\frac{1}{q}\cdot\big(\text{var}(\phi^q)\xi\cdot\xi\big)
\in\Q.
$$ 
Here $\xi$ is a generator of $H_*(M',\p M')$, with $M'$ a connected
component of $M$, 
$\text{var}(\phi^q):H_*(M',\p M')\to H_*(M')$ is the variation
homomorphism of $\phi^q$ and the dot stands for the intersection pairing. 
\smallskip\\
(iii)
An {\bf admissible triple} $(\phi,M,C)$ consists of an orientation
preserving diffeomorphism $\phi:M\to M$ and a $\phi$-invariant finite
union $C\subset M$ of disjoint non-contractible circles, such that the 
following holds.
Let $M'$ be a $\phi$-component of $(M,C)$ and $\chi$ denote
its Euler characteristic.
\smallskip\\
(1) 
If $\chi=0$, then $\phi|\text{cl}(M')$ is an admissible twist map.
\smallskip\\
(2)
If $\chi<0$, then $\phi|\text{cl}(M')$ is a periodic map.
\smallskip\\
By the {\bf period} of a periodic diffeomorphism $\phi$, we mean the
smallest integer $\ell>0$ such that $\phi^\ell=\id$. 
\smallskip\\
(iv) 
Set 
$$
\Tc:=
\{(\chi,d,h;\ell)\in\Z^3\times\Q:d,h>0;\chi\leq0;\chi<0\then\ell\in\N_{>0}\}.
$$
Let $(\phi,M,C)$ be an admissible triple. 
For every $\phi$-component $M'$ of $(M,C)$, set
$t(M'):=(\chi,d,h;\ell)\in\Tc$, where $(\chi,d,h)$ denotes the topological
type of $M'$ and $\ell$ either the period of $\phi|M'$ if $\chi<0$ or the twist
number of $\phi|M'$ if $\chi=0$.
The {\bf characteristic set} of $(\phi,M,C)$ is defined by 
$$
\tc(\phi,M,C):=\big\{t(M'):M'\text{ is a $\phi$-component of }(M,C)\big\}\subset\Tc.
$$
\end{definition}
Recall that an isolated plane curve singularity, or simply 
{\bf isolated singularity}, 
is a germ $[f]$ of holomorphic functions $f:(U,0)\to(\C,0)$, where
$U\subset\C^2$ is a neighborhood of $0$, 
with $(\de f)^{-1}(0)=\{0\}$. The Milnor fiber of $[f]$ is a compact
connected oriented 2-manifold $M$ with boundary. The geometric
monodromy $g$ of $[f]$ is an isotopy class of $\diff_c^+(M)$. 

We now proceed as follows: we explain in \ref{subse:splice}
how to associate a diagram $\Gamma$, called splice diagram, to $[f]$ and 
in \ref{subse:set}, how to associate a set $\tc_\Gamma\subset\Tc$ to such a
diagram $\Gamma$.
The significance of these constructions is expressed by the following
result from \cite{EN}, which we state precisely in Theorem~\ref{thm:EN}:
There exists an admissible triple $(\phi,M,C)$ with characteristic set
$\tc_\Gamma$ and such that $\phi\in\iota_*g$
\footnote{Here, $\iota:\diff^+_c(M)\to\diff^+(M,\p M)$ is the inclusion.}.
Moreover, we have a formula for the twist map components of $\phi$. 
This is used in \ref{subse:proof} to prove Proposition~\ref{prop:monodromy}.

As the title of the monograph \cite{EN} indicates, Eisenbud and Neumann's
point of view is that of link theory. 
Our discussion of splice diagrams, however, is without any
reference to link theory. We simply consider them as a tool for
encoding some algebraic data. This has the advantage that the reader is not 
assumed to be familiar with 3-manifold theory. The downside is that
these diagrams seem to lack intrinsic geometric relevance and 
that it is not clear, how the stated results are actually proven. 
For this we refer to the original literature \cite{EN}, \cite{Ne}. 

\subsection{Puiseux data and splice diagrams}
\label{subse:splice}

Let $[f]$ be an isolated singularity. 
In the following, we define the diagram $\Gamma[f]$ as shown in
\cite[Appendix 1]{EN}. The construction can be divided into 4 steps.  
\begin{step}{1}
The first step uses the so-called Newton method for
solving the equation $f(x,y)=0$ for $y$ in terms of $x$ in a neighborhood of $0$.
We only state the result without going into detail and refer the
reader to \cite[Chapter 8.3]{BK}. 
Recall that a {\bf fractional power series} is a pair $(P,d)$, where
$$
P(x)=\sum_{i=1}^ra_ix^{n_i},
r\in\N\cup\{\infty\},a_i\in\C,a_i\ne0,n_i\in\N_{>0},n_i<n_{i+1},
$$
is a power series converging in a neighborhood of $0$, 
and $d$ is a positive integer which is relatively prime to the set
$\{n_i:i\in\N,i\leq r\}$. Two fractional power series $(P,d),(\tilde{P},\tilde{d})$
are called equivalent, we write $(P,d)\sim (\tilde{P},\tilde{d})$, if
$d=\tilde{d}$ and there exists $\theta\in\C$ such that $\theta^d=1$ and
$\tilde{P}(x)=P(\theta\cdot x)$. 

Now if $[f]$ is an isolated singularity, there is a collection
$\{(P_1,d_1),\dots,(P_\kappa,d_\kappa)\}$ of pairwise non-equivalent
fractional power series, called {\bf Puiseux series}, such that  
$$
f(x,y)=0 \quad\Iff\quad 
\exists\; z\in\C,\exists\; j\in\{1,\dots,\kappa\}:x=z^{d_j},y=P_j(z).
$$
The Puiseux series are uniquely determined up to equivalence;
$\kappa$ is the number of branches of $[f]$.
\end{step}
\begin{step}{2}
Given a collection $\{(P_1,d_1),\dots,(P_\kappa,d_\kappa)\}$ of pairwise
non-equivalent fractional power series, we now define another such
collection $\{(P'_1,d_1),\dots,(P'_\kappa,d_\kappa)\}$ with the property
that each $P'_j$ is a finite series. 
For this we need the following notation. 
Let $\Pi=(P,d)$ be a fractional power series. For $s\in\N$, set 
$$
d_s(\Pi) := \min\{d'\in\N:1\leq i\leq\min(s,r)\then d'\cdot n_i\in d\cdot\N\}.
$$
and 
$$
P^{(s)}(x) := \textstyle\sum_{i=1}^{\min(s,r)} a_ix^{m_i},
\quad
m_i := n_id_s(\Pi)/d.
$$ 
Note that $\Pi^{(s)}:=(P^{(s)},d_s(\Pi))$ is a fractional power series and
that $d_s(P,d)$ is increasing in $s$ and eventually equals $d$. 
For every $j=1,\dots,\kappa$, set $\Pi_j:=(P_j,d_j)$ and define
$$
r_j:=\min\{s\in\N:d_{s}(\Pi_j)=d_j \text{ and }
j\ne j'\then \Pi_j^{(s)}\not\sim\Pi_{j'}^{(s)} \}.
$$
Note that $\Pi_j^{(r_j)}\not\sim\Pi_{i}^{(r_i)}$ if $j\ne i$. 
Now let $\{\Pi_1,\dots,\Pi_\kappa\}$ be the Puiseux series of $[f]$.
We call $\{\Pi_1^{(r_1)},\dots,\Pi_\kappa^{(r_\kappa)}\}$ the 
{\bf Puiseux data} of $[f]$. 
\end{step}
\begin{step}{3} 
Using the Puiseux data of $[f]$, one can now define the diagram
$\tilde{\Gamma}[f]$ from which the diagram $\Gamma[f]$ is obtained in
the next step. For the sake of simplicity, we give the precise
definition of $\tilde{\Gamma}[f]$ only for the cases $\kappa=1,2$. 
Assume that $\kappa=1$ and let $(P,d)$ denote the Puiseux data. 
Define the integers $q_1,\dots,q_r,p_1,\dots,p_r>0$ such that
$\gcd(q_i,p_i)=1$ and 
$$
P(x^{1/d}) =  
x^{\frac{q_1}{p_1}}(a_1+x^{\frac{q_2}{p_1p_2}}
(\dots(a_{r-1}+a_rx^{\frac{q_r}{p_1\cdots p_r}})\dots)).
$$ 
Define the integers $\alpha_1,\dots,\alpha_r$ recursively by 
$$
\alpha_1=q_1,\quad \alpha_{i+1}=p_ip_{i+1}\alpha_i+q_i.
$$
Note that $\gcd(\alpha_i,p_i)=1$ for all $i=1,\dots,r$.
The graph $\tilde{\Gamma}[f]$ is given by 
\medskip
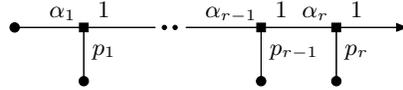
\begin{center}
\begin{picture}(0,0)%
\includegraphics{puiseux1.pstex}%
\end{picture}%
\setlength{\unitlength}{1776sp}%
\begingroup\makeatletter\ifx\SetFigFont\undefined%
\gdef\SetFigFont#1#2#3#4#5{%
  \reset@font\fontsize{#1}{#2pt}%
  \fontfamily{#3}\fontseries{#4}\fontshape{#5}%
  \selectfont}%
\fi\endgroup%
\begin{picture}(5592,1145)(2106,-4266)
\put(2701,-3286){\makebox(0,0)[lb]{\smash{\SetFigFont{5}{6.0}{\familydefault}{\mddefault}{\updefault}{\color[rgb]{0,0,0}\small $\alpha_1$}%
}}}
\put(3376,-3286){\makebox(0,0)[lb]{\smash{\SetFigFont{5}{6.0}{\familydefault}{\mddefault}{\updefault}{\color[rgb]{0,0,0}\small $1$}%
}}}
\put(6901,-3286){\makebox(0,0)[lb]{\smash{\SetFigFont{5}{6.0}{\familydefault}{\mddefault}{\updefault}{\color[rgb]{0,0,0}\small $1$}%
}}}
\put(5776,-3811){\makebox(0,0)[lb]{\smash{\SetFigFont{5}{6.0}{\familydefault}{\mddefault}{\updefault}{\color[rgb]{0,0,0}\small $p_{r-1}$}%
}}}
\put(3301,-3811){\makebox(0,0)[lb]{\smash{\SetFigFont{5}{6.0}{\familydefault}{\mddefault}{\updefault}{\color[rgb]{0,0,0}\small $p_1$}%
}}}
\put(6826,-3811){\makebox(0,0)[lb]{\smash{\SetFigFont{5}{6.0}{\familydefault}{\mddefault}{\updefault}{\color[rgb]{0,0,0}\small $p_r$}%
}}}
\put(6226,-3286){\makebox(0,0)[lb]{\smash{\SetFigFont{5}{6.0}{\familydefault}{\mddefault}{\updefault}{\color[rgb]{0,0,0}\small $\alpha_r$}%
}}}
\put(5851,-3286){\makebox(0,0)[lb]{\smash{\SetFigFont{5}{6.0}{\familydefault}{\mddefault}{\updefault}{\color[rgb]{0,0,0}\small $1$}%
}}}
\put(4876,-3286){\makebox(0,0)[lb]{\smash{\SetFigFont{5}{6.0}{\familydefault}{\mddefault}{\updefault}{\color[rgb]{0,0,0}\small $\alpha_{r-1}$}%
}}}
\end{picture}
\end{center}
\medskip
The coefficients $a_i$ only enter the description of $\tilde{\Gamma}[f]$ if $\kappa>1$. 
Let $\{\Pi,\tilde{\Pi}\}$ be the Puiseux data of $[f]$. To give the
diagram $\tilde{\Gamma}[f]$ one distinguishes three cases. Set 
$$
t := \min\big\{s\geq0:\Pi^{(s+1)}\not\sim\tilde{\Pi}^{(s+1)}\big\}
$$
and let $\tilde{r},\tilde{\alpha}_1,\dots,\tilde{\alpha}_{\tilde{r}}$,
$\tilde{p}_1,\dots,\tilde{p}_{\tilde{r}}$
denote the integers associated to $\tilde{\Pi}$.
\smallskip\\
(i) Assume that $t<r,\tilde{r}$ and $q_{t+1}=\tilde{q}_{t+1},p_{t+1}=\tilde{p}_{t+1}$.
\medskip
\begin{center}
\begin{picture}(0,0)%
\includegraphics{puiseux2.pstex}%
\end{picture}%
\setlength{\unitlength}{1776sp}%
\begingroup\makeatletter\ifx\SetFigFont\undefined%
\gdef\SetFigFont#1#2#3#4#5{%
  \reset@font\fontsize{#1}{#2pt}%
  \fontfamily{#3}\fontseries{#4}\fontshape{#5}%
  \selectfont}%
\fi\endgroup%
\begin{picture}(10842,2420)(1881,-4866)
\put(8401,-3886){\makebox(0,0)[lb]{\smash{\SetFigFont{5}{6.0}{\familydefault}{\mddefault}{\updefault}{\color[rgb]{0,0,0}\small $1$}%
}}}
\put(11251,-3886){\makebox(0,0)[lb]{\smash{\SetFigFont{5}{6.0}{\familydefault}{\mddefault}{\updefault}{\color[rgb]{0,0,0}\small $\tilde{\alpha}_{\tilde{r}}$}%
}}}
\put(9901,-3886){\makebox(0,0)[lb]{\smash{\SetFigFont{5}{6.0}{\familydefault}{\mddefault}{\updefault}{\color[rgb]{0,0,0}\small $\tilde{\alpha}_{\tilde{r}-1}$}%
}}}
\put(8401,-2611){\makebox(0,0)[lb]{\smash{\SetFigFont{5}{6.0}{\familydefault}{\mddefault}{\updefault}{\color[rgb]{0,0,0}\small $1$}%
}}}
\put(11251,-2611){\makebox(0,0)[lb]{\smash{\SetFigFont{5}{6.0}{\familydefault}{\mddefault}{\updefault}{\color[rgb]{0,0,0}\small $\alpha_r$}%
}}}
\put(9901,-2611){\makebox(0,0)[lb]{\smash{\SetFigFont{5}{6.0}{\familydefault}{\mddefault}{\updefault}{\color[rgb]{0,0,0}\small $\alpha_{r-1}$}%
}}}
\put(8251,-4411){\makebox(0,0)[lb]{\smash{\SetFigFont{5}{6.0}{\familydefault}{\mddefault}{\updefault}{\color[rgb]{0,0,0}\small $\tilde{p}_{t+2}$}%
}}}
\put(2476,-3286){\makebox(0,0)[lb]{\smash{\SetFigFont{5}{6.0}{\familydefault}{\mddefault}{\updefault}{\color[rgb]{0,0,0}\small $\alpha_1$}%
}}}
\put(3151,-3286){\makebox(0,0)[lb]{\smash{\SetFigFont{5}{6.0}{\familydefault}{\mddefault}{\updefault}{\color[rgb]{0,0,0}\small $1$}%
}}}
\put(5551,-3286){\makebox(0,0)[lb]{\smash{\SetFigFont{5}{6.0}{\familydefault}{\mddefault}{\updefault}{\color[rgb]{0,0,0}\small $1$}%
}}}
\put(3001,-3811){\makebox(0,0)[lb]{\smash{\SetFigFont{5}{6.0}{\familydefault}{\mddefault}{\updefault}{\color[rgb]{0,0,0}\small $p_1$}%
}}}
\put(5476,-3811){\makebox(0,0)[lb]{\smash{\SetFigFont{5}{6.0}{\familydefault}{\mddefault}{\updefault}{\color[rgb]{0,0,0}\small $p_{t}$}%
}}}
\put(4876,-3286){\makebox(0,0)[lb]{\smash{\SetFigFont{5}{6.0}{\familydefault}{\mddefault}{\updefault}{\color[rgb]{0,0,0}\small $\alpha_{t}$}%
}}}
\put(6001,-3811){\makebox(0,0)[lb]{\smash{\SetFigFont{5}{6.0}{\familydefault}{\mddefault}{\updefault}{\color[rgb]{0,0,0}\small $p_{t+1}$}%
}}}
\put(6931,-3826){\makebox(0,0)[lb]{\smash{\SetFigFont{5}{6.0}{\familydefault}{\mddefault}{\updefault}{\color[rgb]{0,0,0}\small $1$}%
}}}
\put(8251,-3136){\makebox(0,0)[lb]{\smash{\SetFigFont{5}{6.0}{\familydefault}{\mddefault}{\updefault}{\color[rgb]{0,0,0}\small $p_{t+2}$}%
}}}
\put(10726,-3136){\makebox(0,0)[lb]{\smash{\SetFigFont{5}{6.0}{\familydefault}{\mddefault}{\updefault}{\color[rgb]{0,0,0}\small $p_{r-1}$}%
}}}
\put(11776,-3136){\makebox(0,0)[lb]{\smash{\SetFigFont{5}{6.0}{\familydefault}{\mddefault}{\updefault}{\color[rgb]{0,0,0}\small $p_r$}%
}}}
\put(10726,-4411){\makebox(0,0)[lb]{\smash{\SetFigFont{5}{6.0}{\familydefault}{\mddefault}{\updefault}{\color[rgb]{0,0,0}\small $\tilde{p}_{\tilde{r}-1}$}%
}}}
\put(11776,-4411){\makebox(0,0)[lb]{\smash{\SetFigFont{5}{6.0}{\familydefault}{\mddefault}{\updefault}{\color[rgb]{0,0,0}\small $\tilde{p}_{\tilde{r}}$}%
}}}
\put(5866,-3286){\makebox(0,0)[lb]{\smash{\SetFigFont{5}{6.0}{\familydefault}{\mddefault}{\updefault}{\color[rgb]{0,0,0}\small $\alpha_{t+1}$}%
}}}
\put(6931,-3136){\makebox(0,0)[lb]{\smash{\SetFigFont{5}{6.0}{\familydefault}{\mddefault}{\updefault}{\color[rgb]{0,0,0}\small $1$}%
}}}
\put(7351,-4111){\makebox(0,0)[lb]{\smash{\SetFigFont{5}{6.0}{\familydefault}{\mddefault}{\updefault}{\color[rgb]{0,0,0}\small $\tilde{\alpha}_{t+2}$}%
}}}
\put(7321,-2686){\makebox(0,0)[lb]{\smash{\SetFigFont{5}{6.0}{\familydefault}{\mddefault}{\updefault}{\color[rgb]{0,0,0}\small $\alpha_{t+2}$}%
}}}
\put(11851,-2611){\makebox(0,0)[lb]{\smash{\SetFigFont{5}{6.0}{\familydefault}{\mddefault}{\updefault}{\color[rgb]{0,0,0}\small $1$}%
}}}
\put(11851,-3886){\makebox(0,0)[lb]{\smash{\SetFigFont{5}{6.0}{\familydefault}{\mddefault}{\updefault}{\color[rgb]{0,0,0}\small $1$}%
}}}
\put(10801,-2611){\makebox(0,0)[lb]{\smash{\SetFigFont{5}{6.0}{\familydefault}{\mddefault}{\updefault}{\color[rgb]{0,0,0}\small $1$}%
}}}
\put(10801,-3886){\makebox(0,0)[lb]{\smash{\SetFigFont{5}{6.0}{\familydefault}{\mddefault}{\updefault}{\color[rgb]{0,0,0}\small $1$}%
}}}
\end{picture}
\end{center}
\medskip
(ii) Assume that $t<r,\tilde{r}$ and
$\frac{q_{t+1}}{p_{t+1}}<\frac{\tilde{q}_{t+1}}{\tilde{p}_{t+1}}$.
\medskip
\begin{center}
\begin{picture}(0,0)%
\includegraphics{puiseux3.pstex}%
\end{picture}%
\setlength{\unitlength}{1776sp}%
\begingroup\makeatletter\ifx\SetFigFont\undefined%
\gdef\SetFigFont#1#2#3#4#5{%
  \reset@font\fontsize{#1}{#2pt}%
  \fontfamily{#3}\fontseries{#4}\fontshape{#5}%
  \selectfont}%
\fi\endgroup%
\begin{picture}(10842,2420)(1881,-4866)
\put(8401,-3886){\makebox(0,0)[lb]{\smash{\SetFigFont{5}{6.0}{\familydefault}{\mddefault}{\updefault}{\color[rgb]{0,0,0}\small $1$}%
}}}
\put(11251,-3886){\makebox(0,0)[lb]{\smash{\SetFigFont{5}{6.0}{\familydefault}{\mddefault}{\updefault}{\color[rgb]{0,0,0}\small $\tilde{\alpha}_{\tilde{r}}$}%
}}}
\put(8401,-2611){\makebox(0,0)[lb]{\smash{\SetFigFont{5}{6.0}{\familydefault}{\mddefault}{\updefault}{\color[rgb]{0,0,0}\small $1$}%
}}}
\put(11251,-2611){\makebox(0,0)[lb]{\smash{\SetFigFont{5}{6.0}{\familydefault}{\mddefault}{\updefault}{\color[rgb]{0,0,0}\small $\alpha_r$}%
}}}
\put(9901,-2611){\makebox(0,0)[lb]{\smash{\SetFigFont{5}{6.0}{\familydefault}{\mddefault}{\updefault}{\color[rgb]{0,0,0}\small $\alpha_{r-1}$}%
}}}
\put(8251,-4411){\makebox(0,0)[lb]{\smash{\SetFigFont{5}{6.0}{\familydefault}{\mddefault}{\updefault}{\color[rgb]{0,0,0}\small $\tilde{p}_{t+1}$}%
}}}
\put(2476,-3286){\makebox(0,0)[lb]{\smash{\SetFigFont{5}{6.0}{\familydefault}{\mddefault}{\updefault}{\color[rgb]{0,0,0}\small $\alpha_1$}%
}}}
\put(3151,-3286){\makebox(0,0)[lb]{\smash{\SetFigFont{5}{6.0}{\familydefault}{\mddefault}{\updefault}{\color[rgb]{0,0,0}\small $1$}%
}}}
\put(5551,-3286){\makebox(0,0)[lb]{\smash{\SetFigFont{5}{6.0}{\familydefault}{\mddefault}{\updefault}{\color[rgb]{0,0,0}\small $1$}%
}}}
\put(3001,-3811){\makebox(0,0)[lb]{\smash{\SetFigFont{5}{6.0}{\familydefault}{\mddefault}{\updefault}{\color[rgb]{0,0,0}\small $p_1$}%
}}}
\put(5476,-3811){\makebox(0,0)[lb]{\smash{\SetFigFont{5}{6.0}{\familydefault}{\mddefault}{\updefault}{\color[rgb]{0,0,0}\small $p_{t}$}%
}}}
\put(4876,-3286){\makebox(0,0)[lb]{\smash{\SetFigFont{5}{6.0}{\familydefault}{\mddefault}{\updefault}{\color[rgb]{0,0,0}\small $\alpha_{t}$}%
}}}
\put(8251,-3136){\makebox(0,0)[lb]{\smash{\SetFigFont{5}{6.0}{\familydefault}{\mddefault}{\updefault}{\color[rgb]{0,0,0}\small $p_{t+2}$}%
}}}
\put(10726,-3136){\makebox(0,0)[lb]{\smash{\SetFigFont{5}{6.0}{\familydefault}{\mddefault}{\updefault}{\color[rgb]{0,0,0}\small $p_{r-1}$}%
}}}
\put(11776,-3136){\makebox(0,0)[lb]{\smash{\SetFigFont{5}{6.0}{\familydefault}{\mddefault}{\updefault}{\color[rgb]{0,0,0}\small $p_r$}%
}}}
\put(10726,-4411){\makebox(0,0)[lb]{\smash{\SetFigFont{5}{6.0}{\familydefault}{\mddefault}{\updefault}{\color[rgb]{0,0,0}\small $\tilde{p}_{\tilde{r}-1}$}%
}}}
\put(11776,-4411){\makebox(0,0)[lb]{\smash{\SetFigFont{5}{6.0}{\familydefault}{\mddefault}{\updefault}{\color[rgb]{0,0,0}\small $\tilde{p}_{\tilde{r}}$}%
}}}
\put(5866,-3286){\makebox(0,0)[lb]{\smash{\SetFigFont{5}{6.0}{\familydefault}{\mddefault}{\updefault}{\color[rgb]{0,0,0}\small $\alpha_{t+1}$}%
}}}
\put(6931,-3136){\makebox(0,0)[lb]{\smash{\SetFigFont{5}{6.0}{\familydefault}{\mddefault}{\updefault}{\color[rgb]{0,0,0}\small $1$}%
}}}
\put(7351,-4111){\makebox(0,0)[lb]{\smash{\SetFigFont{5}{6.0}{\familydefault}{\mddefault}{\updefault}{\color[rgb]{0,0,0}\small $\tilde{\alpha}_{t+1}$}%
}}}
\put(7321,-2686){\makebox(0,0)[lb]{\smash{\SetFigFont{5}{6.0}{\familydefault}{\mddefault}{\updefault}{\color[rgb]{0,0,0}\small $\alpha_{t+2}$}%
}}}
\put(6676,-3811){\makebox(0,0)[lb]{\smash{\SetFigFont{5}{6.0}{\familydefault}{\mddefault}{\updefault}{\color[rgb]{0,0,0}\small $p_{t+1}$}%
}}}
\put(11851,-2611){\makebox(0,0)[lb]{\smash{\SetFigFont{5}{6.0}{\familydefault}{\mddefault}{\updefault}{\color[rgb]{0,0,0}\small $1$}%
}}}
\put(11851,-3886){\makebox(0,0)[lb]{\smash{\SetFigFont{5}{6.0}{\familydefault}{\mddefault}{\updefault}{\color[rgb]{0,0,0}\small $1$}%
}}}
\put(9901,-3886){\makebox(0,0)[lb]{\smash{\SetFigFont{5}{6.0}{\familydefault}{\mddefault}{\updefault}{\color[rgb]{0,0,0}\small $\tilde{\alpha}_{\tilde{r}-1}$}%
}}}
\put(10801,-2611){\makebox(0,0)[lb]{\smash{\SetFigFont{5}{6.0}{\familydefault}{\mddefault}{\updefault}{\color[rgb]{0,0,0}\small $1$}%
}}}
\put(10801,-3886){\makebox(0,0)[lb]{\smash{\SetFigFont{5}{6.0}{\familydefault}{\mddefault}{\updefault}{\color[rgb]{0,0,0}\small $1$}%
}}}
\end{picture}
\end{center}
\medskip
(iii) Assume that $t=\tilde{r}<r$. In this case, the diagram is obtained from the
diagram in case (ii) by terminating the edge with weight $p_{t+1}$ by an
arrowhead. 
\smallskip\\
By interchanging $\Pi$ and $\tilde{\Pi}$ if necessary, (i--iii) define
$\tilde{\Gamma}[f]$ in the case $\kappa=2$. 
The general case is obtained by induction on $\kappa$. The induction step
involves operations of the kind (i--iii). 
Instead of giving the precise definition, we give a list of
properties of the diagram $\tilde{\Gamma}=\tilde{\Gamma}[f]$ for an
arbitrary singularity $[f]$. In the case $\kappa=1,2$, these properties are
easily verified from the definitions above.  
\smallskip\\
(A1)
$\tilde{\Gamma}$ has the structure of a weighted tree. All weights are
positive integers. 
\smallskip\\
(A2)
$\tilde{\Gamma}$ has three kinds of vertices: arrowhead, knob and
box vertices.  
The number of arrowhead vertices equals the number of branches of $[f]$. 
The arrowhead and knob vertices have $1$ incoming edge, the boxed ones at
least $3$. A box vertex has at most $2$ neighboring knob vertices. 
\smallskip \\
(A3)
An edge of $\tilde{\Gamma}$ carries a weight at each ending box vertex.
The edge-weights at a box vertex are pairwise relatively prime. 
\smallskip\\
(A4)
Let $b$ denote a box vertex of $\tilde{\Gamma}$. Let $E$ be the set of edges
which connect $b$ to its neighboring box vertices. 
The graph $\tilde{\Gamma}\setminus E$ has $\# E+1$ connected components. 
There is at most one component which does neither contain $b$ nor any arrowhead
vertex. 
Assume that there is such a component. 
Let $e$ be the edge which connects $b$ to that component and let $b'$ denote
the other vertex of $e$. Then $e$ has weight $1$ at $b'$. 
\smallskip\\
(A5)
Let $b,b'$ denote connected box vertices of $\tilde{\Gamma}$.  
Let $a_1,\dots,a_k$ be the weights at $b$ of the edges connecting $b$
and its neighboring vertices. Similarly, let $a'_1,\dots,a'_{k'}$ be the
weights at $b'$. Assume that $a_1,a'_1$ are the weights of the edge connecting
$b$ and $b'$. Then 
$$
a_1a'_1-a_2\cdots a_ka'_2\cdots a'_{k'}>0.
$$
\end{step}
\begin{step}{4}
The diagram $\Gamma=\Gamma[f]$ is obtained
from $\tilde{\Gamma}=\tilde{\Gamma}[f]$ by the following algorithm.
Set $\Gamma_0:=\tilde{\Gamma}$.
{\sc Step}: 
Let $e$ be an edge of $\Gamma_i$ connecting a box vertex $b$ and
a knob vertex $v$ and having weight $1$.
If $e$ does not exist, set $\Gamma:=\Gamma_i$ and {\sc Stop}.
Otherwise, let $k\geq3$ denote the number of incoming edges at
$b$ and $n\leq2$ the number of neighboring box vertices of $b$.
If $k=3$ and $n=0$, set $\Gamma:=\Gamma_*:=\includegraphics{gamma1.pstex}$ 
and {\sc Stop}.

Otherwise, apply the operation $\Gamma_i\to\Gamma_{i+1}$ defined as follows
\medskip\\
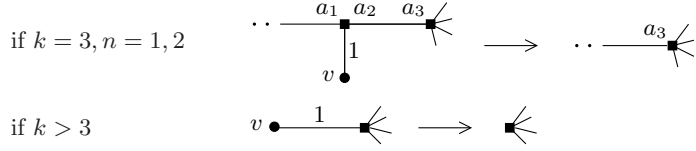
\begin{minipage}{3.5cm}
\small $\quad$ if $k=3,n=1,2$
\end{minipage}
\begin{minipage}{7cm}
\begin{picture}(0,0)%
\includegraphics{operation12n.pstex}%
\end{picture}%
\setlength{\unitlength}{1776sp}%
\begingroup\makeatletter\ifx\SetFigFont\undefined%
\gdef\SetFigFont#1#2#3#4#5{%
  \reset@font\fontsize{#1}{#2pt}%
  \fontfamily{#3}\fontseries{#4}\fontshape{#5}%
  \selectfont}%
\fi\endgroup%
\begin{picture}(6187,1175)(3511,-4191)
\put(9001,-3511){\makebox(0,0)[lb]{\smash{\SetFigFont{5}{6.0}{\familydefault}{\mddefault}{\updefault}{\color[rgb]{0,0,0}\small $a_3$}%
}}}
\put(4876,-3826){\makebox(0,0)[lb]{\smash{\SetFigFont{5}{6.0}{\familydefault}{\mddefault}{\updefault}{\color[rgb]{0,0,0}\small $1$}%
}}}
\put(4426,-3211){\makebox(0,0)[lb]{\smash{\SetFigFont{5}{6.0}{\familydefault}{\mddefault}{\updefault}{\color[rgb]{0,0,0}\small $a_1$}%
}}}
\put(4951,-3211){\makebox(0,0)[lb]{\smash{\SetFigFont{5}{6.0}{\familydefault}{\mddefault}{\updefault}{\color[rgb]{0,0,0}\small $a_2$}%
}}}
\put(5626,-3211){\makebox(0,0)[lb]{\smash{\SetFigFont{5}{6.0}{\familydefault}{\mddefault}{\updefault}{\color[rgb]{0,0,0}\small $a_3$}%
}}}
\put(4516,-4156){\makebox(0,0)[lb]{\smash{\SetFigFont{5}{6.0}{\familydefault}{\mddefault}{\updefault}{\color[rgb]{0,0,0}\small\bf $v$}%
}}}
\end{picture}
\end{minipage}
\medskip\\
\begin{minipage}{3.5cm}
\small $\quad$ if $k>3$
\end{minipage}
\begin{minipage}{5cm}
\begin{picture}(0,0)%
\includegraphics{operation3.pstex}%
\end{picture}%
\setlength{\unitlength}{1776sp}%
\begingroup\makeatletter\ifx\SetFigFont\undefined%
\gdef\SetFigFont#1#2#3#4#5{%
  \reset@font\fontsize{#1}{#2pt}%
  \fontfamily{#3}\fontseries{#4}\fontshape{#5}%
  \selectfont}%
\fi\endgroup%
\begin{picture}(3997,569)(4876,-3908)
\put(4876,-3661){\makebox(0,0)[lb]{\smash{\SetFigFont{5}{6.0}{\familydefault}{\mddefault}{\updefault}{\color[rgb]{0,0,0}\small $v$}%
}}}
\put(5761,-3541){\makebox(0,0)[lb]{\smash{\SetFigFont{5}{6.0}{\familydefault}{\mddefault}{\updefault}{\color[rgb]{0,0,0}\small $1$}%
}}}
\end{picture}
\end{minipage}
\medskip\\
Finally repeat the {\sc Step}. 

Among all the possible diagrams $\Gamma$ which are obtained by this algorithm, 
$\Gamma_*$ is exceptional in the sense that it does not have box vertices.  
It is important to note, that in all other cases the properties
(A1--5) still hold if $\tilde{\Gamma}$ is replaced by $\Gamma$. 
Additionally, we have 
\smallskip\\
(A6)
Every edge of $\Gamma$ connecting a box and a knob vertex has weight $>1$.   
\end{step}
\smallskip\\
We end the discussion of splice diagrams with the remark that the diagram 
associated to the quadratic singularity $(x,y)\mapsto x^2+y^2$ is $\Gamma_*$.

\subsection{Characteristic set}
\label{subse:set}

Let $\Gamma$ be the splice diagram of an isolated singularity.
In this section, we define the set $\tc_\Gamma\subset\Tc$. 
In the exceptional case $\Gamma=\Gamma_*$, this is simply 
the set $\{(0,1,2;1)\}$, which is the characteristic set of
a positive Dehn twist. Assume from now on that $\Gamma\ne\Gamma_*$. 
Denote by $\Ac$ respectively $\Bc$ the set of arrowhead respectively box
vertices of $\Gamma$.  
Moreover, denote by $\Ec$ the set of
edges of $\Gamma$ which connect $\Bc$ to $\Ac\cup\Bc$. 

In the following, we define 
{\bf (i)} for each $b\in\Bc$ an element $t_b\in\Tc$ and {\bf (ii)} for each
$e\in\Ec$ an element $t_e\in\Tc$. We then set 
$$
\tc_\Gamma:=\{t_x:x\in\Bc\cup\Ec\}.
$$
Finally we define {\bf (iii)} for each $e\in\Ec$ an admissible twist map
$\phi_e$. 

Let $\Vc$ denote the set of ordered pairs of connected vertices of
$\Gamma$.
We start by introducing the function $m:\Vc\to\N$. 
Let $(v,v')\in\Vc$ and let $e$ be the edge connecting $v$ and $v'$. 
The graph $\Gamma\setminus\{e\}$ has two components.
Denote by $\Gamma'$ that component which contains $v'$. 
Let $\Ac'$ denote the set of arrowhead vertices of $\Gamma'$.
For each $a\in\Ac'$, there is a unique path in $\Gamma'$, denoted by $\gamma_a$,
which connects $v'$ and $a$.  
Define $\sigma_a>0$ to be the product of all edge-weights adjacent to
$\gamma_a$, but not on $\gamma_a$. For examples see \cite[page 84]{EN}. 
If $\gamma_a$ is the constant path, set $\sigma_a:=1$. Define
$$
m(v,v'):= \sum_{a\in\Ac'}\sigma_a.
$$
Note that $m(v,v')=0$ if and only if $\Ac'=\emptyset$.
Together with (A4), this has the following consequences:  
\smallskip\\
(B1)
If $m(v,v')=0$, then $m(v',v)\ne0$.
\smallskip\\
(B2)
If $b,b'\in\Bc$ and $m(b,b')=0$, then the edge connecting $b$ and $b'$ has
weight $1$ at $b'$. 
\smallskip\\
(B3)
For each $b\in\Bc$, there exists at most one neighboring vertex
$b'\in\Bc$ such that $m(b,b')=0$.  
\smallskip\\
(B4)
If $b\in\Bc$ and $a\in\Ac$, then $m(b,a)=1$. 
\smallskip\\
{\bf (i)}
Let $b\in\Bc$ and denote by $\{v_1,\dots,v_n,v_{n+1},\dots,v_k\},1\leq n\leq k,$ 
the set of vertices which are connected to $b$, ordered such that
$v_i\in\Ac\cup\Bc$ if and only if $i\leq n$. 
Note that $k\geq 3$ and that $k-n\in\{0,1,2\}$. 
Further denote by
$a_i>0,i=1,\dots,k,$ the weight at $b$ of the edge connecting $b$ and
$v_i$. Define the numbers
$$
d_b:=\gcd\big(m(b,v_1),\dots,m(b,v_n)\big),\quad
h_b:=\sum_{i=1}^n\gcd\big(m(b,v_i),m(v_i,b)\big)
$$
and
\begin{equation}\label{eq:b}
\ell_b:=\sum_{i=1}^nm(b,v_i)\cdot a_1\cdots\widehat{a_i}\cdots a_k,\quad
\chi_b:=\ell_b\cdot\big(2-k+\sum_{i=n+1}^k\frac{1}{a_i}\big),
\end{equation}
where the hat means that the underlying factor is omitted. 
Note that (B3), (B4) respectively (B1) imply that $d_b$
respectively $h_b$ is well defined and therefore $>0$.  
Similarly, it follows from (B3), (B4) that the integer $\ell_b$ is $>0$. 
We furthermore claim, that $\chi_b<0$. This is because $k\geq3$ and
$a_{n+1},\dots,a_k>1$ are pairwise relatively prime.
Hence, we can set 
$$
t_b:=(\chi_b,d_b,h_b;\ell_b).
$$ 
Finally note, that from the definition of $m$, it follows that for each
$i=1,\dots,n$,
\begin{equation}\label{eq:ellb}
\ell_b = 
m(b,v_i)\cdot a_1\cdots\widehat{a_i}\cdots a_k + m(v_i,b)\cdot a_i.
\end{equation}
{\bf (ii)}
Let $e\in\Ec$ be an edge which connects the vertices $b,b'\in\Bc$. 
Let $a_1,\dots,a_k$ be the weights at $b$ of the edges connecting $b$
and its neighboring vertices. Similarly, let $a'_1,\dots,a'_{k'}$ be the
weights at $b'$. Assume that $a_1,a'_1$ are the weights of
$e$. Define the numbers
\begin{equation}\label{eq:deltae}
d_e:=\gcd\big(m(b,b'),m(b',b)\big),\quad
\Delta_e:=a_1a'_1-a_2\cdots a_ka'_2\cdots a'_{k'}
\end{equation}
and 
$$
\ell_e:=\frac{d_e\cdot\Delta_e}{\ell_b\cdot\ell_{b'}}.
$$
By (B1), $d_e$ is well defined and therefore $>0$. 
Set 
$$
t_e:=(0,d_e,2d_e;\ell_e)\in\Tc.
$$ 
Now let $e\in\Ec$ be an edge which connects the vertex $b\in\Bc$ to an
arrowhead vertex. In this case set
$$
t_e := (0,d_e,2d_e;\ell_e) := (0,1,2;1/\ell_b)\in\Tc.
$$
{\bf (iii)}
Let $e$ be an edge which connects $b\in\Bc$ and $b'\in\Ac\cup\Bc$. 
By (B1),(B4) we can assume without loss of generality that $m(b,b')\ne0$. 
Now choose integers $n,n'$ 
with $m(b,b')\cdot n'+m(b',b)\cdot n=d_e$. 
Denote by $a$ the weight of $e$ at $b$ and set $m:=m(b,b')$. 
Define $\phi_e$ to be the admissible twist map which cyclically permutes $d_e$
annuli and such that 
\begin{equation}\label{eq:e}
\phi_e^{d_e}(q,p) =  
\Big(q,p-q\cdot d_e\ell_e-\frac{d_e}{m}\big(n-\frac{d_e\cdot a}{\ell_b}\big)\Big),
\end{equation}
for all $(q,p)\in[0,1]\times S^1$. 
\subsection{Proof of Proposition~\ref{prop:monodromy}}
\label{subse:proof}

We first state the precise results that are needed for the proof. 
\begin{theorem}[Eisenbud, Neumann]\label{thm:EN}
Let $[f]$ be an isolated plane curve singularity with Milnor fiber $M$ 
and geometric monodromy $g$. 
There exists an admissible triple $(\phi,M,C)$ such that $\phi\in\iota_*g$ and 
$\tc(\phi,M,C)=\tc_{\Gamma[f]}$.
Moreover, the twist map components of $\phi$ are given by the
model~\eqref{eq:e}.  
\end{theorem}
\begin{remarks}
(i)
This theorem is a summary of results from Sections~9, 10, 11 and 13
of \cite{EN}. 
The periodic components of the monodromy are described in 
Lemma~11.4 and the twist map components in Theorems~13.1, 13.5. 
The positive twist property (A5) is contained in Theorem~9.4. 
We remark that our notation differs at some points from that of
\cite{EN}. 
\smallskip\\
(ii) Eisenbud, Neumann give a detailed description of the periodic
components of the monodromy in terms of cyclic branched coverings in
\cite[Lemma~11.4]{EN}. 
\smallskip\\
(iii)
The graph $\Gamma[f]$ contains more information than the characteristic set
$\tc_{\Gamma[f]}$. It also shows how the $\phi$-components are pieced together
to give the Milnor fiber. 
\end{remarks}
The second main ingredient for our proof is the following result of \cite{AC2}
and \cite{L}. 
We would like to point out that this result holds in much greater generality
than we use it here,
namely for holomorphic hypersurface singularities, isolated or non-isolated,
in any dimension.
\begin{theorem}[A'Campo, L{\^e}]\label{thm:AL}
Let $g$ be the geometric monodromy of an isolated plane curve singularity.  
Then $\Lambda(g)=0$, where $\Lambda$ denotes the Lefschetz number. 
\end{theorem}
\begin{proof}[Proof of Proposition \ref{prop:monodromy}]
Let $[f]$ be an isolated plane curve singularity. 
Let $M$ denote the Milnor fiber and $g$ the geometric monodromy of $[f]$. 
By Theorem~\ref{thm:EN}, 
there exists a representative $\phi\in\iota_*g$ and a finite $\phi$-invariant
union of circles $C\subset M$, such that if $M'$ is a $\phi$-component of
$M\setminus C$, then either $M'$ has negative Euler characteristic and $\phi|M'$ is
periodic, or $M'$ is a union of annuli and $\phi|M'$ is an admissible twist map. 
If $(\chi,d,h)$ denotes the topological type of $M'$ and $\ell$ the
order/twist number of $\phi|M'$, then $(d,\chi,h;\ell)\in\tc_{\Gamma[f]}$. 
Moreover, if $M'$ is a union of annuli, then \eqref{eq:e} is a model
for $\phi|M'$. 
 
The strategy of the proof is now as follows.
We will prove below that 
\begin{claim}{1}
If $\chi=0$, then $\fix(\phi)\cap\text{int}(M')=\emptyset$.
\end{claim}
\begin{claim}{2}
If $\chi<0$, then $\ell>1$. 
\end{claim}
\begin{claim}{3}
$\phi$ only has positive twists.  
\end{claim}
\smallskip\\
Claim~1 implies that $\phi$ is a diffeomorphism of finite type. 
From Proposition~\ref{prop:fclass} about the fixed point classes of a
diffeomorphism of finite type and claim~2, 
it follows that $\fix(\phi)\cap\text{int}(M)$ is a discrete set of fixed
points with fixed point index~$1$. The Lefschetz fixed point theorem therefore
implies that 
$$
\Lambda(\phi)=\#\big(\fix(\phi)\cap\text{int}(M)\big).
$$
Note that $\p M$ has fixed point index $0$. 
From Theorem~\ref{thm:AL}, it hence follows that
$\fix(\phi)\cap\text{int}(M)=\emptyset$.
Together with claim~3, this proves Proposition~\ref{prop:monodromy}, up to 
claim~1,2 and 3.

Note that if $\Gamma[f]=\Gamma_*$, then $M$ is an annulus and $\phi$
is a positive Dehn twist. Claim~1,2 and 3 obviously hold in this case
and we assume from now on that $\Gamma[f]\ne\Gamma_*$.
We begin by proving claim~2. 
Recall that $\tc_{\Gamma[f]}=\{t_x:x\in\Bc\cup\Ec\}$ and that if 
$\chi<0$, there exists $b\in\Bc$ such that
$(\chi,d,h;\ell)=(\chi_b,d_b,h_b;\ell_b)$. 
Consider equation~\eqref{eq:b} and 
assume that $\ell_b=1$. In this case, only one summand of $\ell_b$ is
non-zero, which, by (B3), is only possible if $k\leq2$.   
Since $n\geq3$, however, it follows from (A6) that $a_n>1$ and hence that
$1=\ell_b\geq a_n>1$, a contradiction.

To proof claim~1, let $e\in\Ec$ be such that 
$(\chi,d,h;\ell)=(0,d_e,2d_e;\ell_e)$. We can assume that $d_e=1$, otherwise
the claim is obviously true. Hence $M'$ is an annulus and $\phi|M'=\phi_e$.
Consider equation~\eqref{eq:e} and assume that $\phi_e(q,p)=(q,p)$ for some 
$(q,p)\in(0,1)\times S^1$. 
This is only possible if 
$$
-q\cdot\ell_e-\frac{1}{m}\big(n-\frac{a}{\ell_b}\big)
\in\Z,
$$
where we use the same notation as in \eqref{eq:e}.
This in turn implies that 
\begin{equation}\label{eq:q}
a-qm\ell_b\ell_e\in\ell_b\Z.
\end{equation}
Assume for the moment that $0<m\ell_b\ell_e\leq a$.
Since $0<q<1$, it follows from \eqref{eq:q} that $0<\ell_b<a$.
To show that this is a contradiction, recall from \eqref{eq:ellb} that 
$\ell_b=ma_2\cdots a_r+m'a$, where $m':=m(b',b)$. 
If $m'\ne0$, it follows that $\ell_b\geq a$ and
hence $m'=0$. 
By (B2) however, this implies that $a=1>\ell_b$, which is a contradiction.

It remains to prove that $0<m\ell_b\ell_e\leq a$.
First note, that $m\ell_b\ell_e>0$ iff $\ell_e>0$. 
If $b'\in\Ac$, then $\ell_b>0$ by definition. 
If $b\in\Bc$, then (A5) is exactly the statement that $\ell_e>0$.
This in fact proves claim~3. 
To prove that $m\ell_b\ell_e\leq a$, first assume that $b'\in\Ac$. 
Then $m=1,a=1,\ell_b\ell_e=1$ and we are finished. 
If $b'\in\Bc$, then it follows from \eqref{eq:ellb} and \eqref{eq:deltae} that
$$
\ell_{b'} \geq ma',
\quad
\Delta_e \leq aa',
$$
where $a'$ denotes the weight of $e$ at $b'$. 
This implies that 
$$
\ell_e\leq\frac{a}{\ell_bm},
$$ 
which proves the required inequality. This ends the proof of the proposition.  
\end{proof}

\section{Floer homology on surfaces with boundary}\label{ap:open}

Let $M$ be a compact connected oriented 2-manifold with boundary 
$\p M\ne\emptyset$. Recall that $\diff_c^+(M)$ denotes the group of
orientation preserving diffeomorphisms which are the identity near $\p M$. 
This appendix addresses Floer homology theory for elements of
$\diff_c^+(M)$.
In higher dimensions, this is known as Floer homology theory for exact
symplectomorphisms of exact symplectic manifolds with contact type boundary
and was used in \cite[Section~4]{S5}, see also \cite{CFH}. 
Due to the dimensional restriction, the theory exhibits auxiliary structure,
namely isotopy invariance, as in the closed case.  
The central notion around this issue is that of monotonicity.
We start by defining monotonicity and show that it has naturality, isotopy and
inclusion properties similar to the ones discussed on
page~\pageref{page:natur} for the closed case. 

Let $\omega$ be an area form on $M$ and denote by $\symp_c(M,\omega)$ the group of
$\omega$-preserving diffeomorphisms which are the identity near $\p M$. 
If $\phi\in\symp_c(M,\omega)$, $\omega$ induces a
closed 2-form $\omega_\phi$ on the mapping torus $T_\phi$. 
\begin{definition}\label{def:open}
$\phi\in\symp_c(M,\omega)$ is called {\bf monotone}, if $[\omega_\phi]=0$ in
$H^2(T_\phi;\R)$. $\symp_c^m(\Sigma,\omega)$ denotes the set of
monotone symplectomorphisms. 
\end{definition}
As in the closed case it is useful to look at the short exact sequence
$$
0 \To  \frac{H^1(M;\R)}{\im(\id-\phi^*)}
\stackrel{\delta}{\To} H^2(\tf;\R)
\stackrel{\iota^*}{\To} H^2(M;\R)=0
\To 0
$$
and define the class $m(\phi)\in H^1(M;\R)/\im(\id-\phi^*)$ satisfying
$\delta m(\phi)=[\omega_\phi]$. 
The naturality, isotopy and inclusion properties discussed on
page~\pageref{page:natur} in the closed case carry over word by word to the current
situation with the addition of a subscript $c$ to all diffeomorphism groups. 
For the first two properties this is straight forward to check. 
The inclusion property needs separate consideration. Recall:
\smallskip\\
(Inclusion)
The inclusion $\symp_c^m(M,\omega)\inclusion\diff_c^+(M)$ is a homotopy
equivalence.
In particular, every connected component of $\symp_c^m(\Sigma,\omega)$ is
contractible. 
\smallskip\\
The proof is analogous to the closed case and uses the following three facts. 
Firstly, the inclusion
$\symp_c(M,\omega)\inclusion\diff_c^+(M)$ is a homotopy equivalence. 
This follows from an extension of Moser's theorem, see~\cite[Exercise~3.18]{MS}. 
Secondly, every connected component of $\diff_c^+(M)$ is contractible.
If the genus of $M$ is $\ne0$, this is shown using the Earl-Eells
Theorem~\cite{EE}. If $g=0$, it follows from the corresponding result for the
disk, which is due to Smale~\cite[Theorem~B]{Sm1}. Thirdly, we have 
\begin{lemma}
If $\phi\in\symp_c(M,\omega)$, there exists 
a closed 1-form $\theta\in m(\phi)$ such that
$\supp(\theta)\subset\mathrm{int}(M)$. The flow 
$(\psi_t)_{t\in\R}$ of the vector field $X$ which is uniquely defined by 
$\omega(X,\cdot)=-\theta$, satisfies $\psi_t\in\symp_c(M,\omega)$ and 
$m(\phi\circ\psi_1)=0$.
\end{lemma}
\begin{proof}
The second part of the statement follows immediately from the first one and
the isotopy property: 
$$
m(\phi\circ\psi_1)=m(\phi)+[\flux(\psi_t)]
$$  
in $H^1(M;\R)/\im(\id-\phi^*)$.
To prove the first statement, let $\beta\in m(\phi)$ be a
closed 1-form. 
Let $S_1,\dots,S_n$ denote the connected components of $\p M$ and
choose for each $i$ a collar neighborhood $N_i\subset M$ of $S_i$ and a closed 1-form
$\theta_i$ on $N_i$ such that $\langle[\theta_i],[S_i]\rangle=1$. 
There exist $f:M\to\R$ smooth and $t_1,\dots,t_n\in\R$ such that 
\begin{equation}\label{eq:open1}
(\beta+\de f)|N_i = t_i\cdot\theta_i,
\quad \forall\, i=1,\dots,n.
\end{equation}
We claim that $\theta:=\beta+\de f$ is the required 1-form. 
Recall from page~\pageref{eq:cohomology}, that $\delta:H^1(M;\R)\to
H^2(\tf;\R)$ is given by 
$\delta[\theta]=[\rho\cdot\theta\wedge\de t]$, with $\rho:\I\to\R$ a smooth
function vanishing near 0 and 1, and satisfying $\int_0^1\!\rho\,\de t=1$.
Furthermore note, that $S^1\times S_i\subset T_\phi$ is an embedded 2-torus
for each $i=1,\dots,n$. 
Using \eqref{eq:open1}, it follows that 
\begin{eqnarray*}
\big\langle[\omega_\phi],[S^1\times S_i]\big\rangle 
&=& \big\langle\delta[\alpha],[S^1\times S_i]\big\rangle \\
&=& -\big\langle[\rho\cdot\de t],[S^1]\big\rangle\cdot
\big\langle[\alpha],[S_i]\big\rangle \\
&=& -t_i\cdot\big\langle[\theta_i],[S_i]\big\rangle
\;=\; -t_i.
\end{eqnarray*}
On the other hand, 
$$
\langle[\omega_\phi],[S^1\times S_i]\rangle=0, 
$$
since
$\omega_\phi$ has no $\de t$-component. Hence, $t_i=0$ for all $i=1,\dots,n$,
which proves the claim. 
\end{proof}
Recall that in the closed case monotonicity (i) guarantees compactness of the
space of Floer connecting orbits and (ii) is used to prove invariance. 
The same holds in the current situation. 
\smallskip\\
(Floer homology)
To every $\phi\in\symp_c^m(\Sigma,\omega)$ symplectic Floer homology
theory assigns a pair of $\Z_2$-graded vector spaces $\HF(\phi,\pm)$
over $\Z_2$, with multiplicative structures
$$
H^*(M;\Z_2)\otimes\HF(\phi,\pm)\To\HF(\phi,\pm).
$$
Floer homology is natural in the sense that $\HF(\phi,\pm)$ and
$\HF(\psi^*\phi,\pm)$ are naturally isomorphic as modules over
$H^*(M;\Z_2)$, for all $\psi\in\diff_c^+(M)$. 
\smallskip\\
(Invariance) 
If $\phi,\phi'\in\symp_c^m(\Sigma,\omega)$ are isotopic, then
$\HF(\phi,\pm)$ and $\HF(\phi',\pm)$ are naturally isomorphic as modules over
$H^*(M;\Z_2)$.
\smallskip\\ 
The appearance of the sign in the Floer homology corresponds to two ways of
perturbing $\phi\in\symp_c^m(\Sigma,\omega)$ near $\p M$. 
To be more precise, let $\jmath:\bigcup(-\e,0]\times S^1\to M$ be a collar
neighborhood of $\p M$, such that $\jmath^*\omega=\de q\wedge\de p$ with
$(q,p)\in(-\e,0]\times S^1$. Choose $H:M\to\R$ with support
near $\p M$ and such that $\jmath^*H(q,p)=-q$. Let $\psi_t$ denote the Hamiltonian
flow generated by $H$, choose $0<\delta<1$ and set 
$$
\phi_+ := \phi\circ\psi_\delta,
\quad
\phi_- := \phi\circ\psi_{-\delta}.
$$
The definition of the Floer complex for $\phi_\pm$ is along the
same line as that in the closed case \cite{S2}, with the usual
modifications that are needed in the presence of a contact type boundary.

The modifications include a condition on the path $J=(J_t)_{t\in\R}$ of
$\omega$-compatible complex structures that is used to define the Floer
connecting orbits; namely that $\jmath^*J_t$ is the standard complex
structure on $\bigcup(-\e,0]\times S^1$, for all $t\in\R$.
We briefly recall the use of this condition. 
Assume without loss of generality that $\phi|\im\jmath=\id$. 
Now let $u:\R^2\to M$ be a smooth map satisfying
\begin{equation}\label{eq:open2}
\left\{\begin{array}{l}
u(s,t) = \phi_+(u(s,t+1)), \\
\p_s u + J_t(u)\p_t u = 0, \\
\lim_{s\to\pm\infty}u(s,t) \in\fix(\phi_+).
\end{array}\right.
\end{equation}
We claim that $\im u\subset M\setminus\im\jmath$.  
Assume by contradiction that $u^{-1}(\im\jmath)$ is non-empty and let
$u_q:u^{-1}(\im\jmath)\to\R$ denote the $q$-component of $\jmath^{-1}\circ u$.
By construction, $u_q$ is smooth and not locally constant. 
Using the first and third equation in \eqref{eq:open2}, one can now
show that $u_q$ has a global maximum.
From the second equation in \eqref{eq:open2} together with the
above assumption on $J_t$, it furthermore follows that $u_q$ is a
harmonic function.
This contradicts the maximum principle and hence
proves the claim, which assures that the Floer connecting orbits are
contained in a compact subset of $\text{int}(M)$.

We close this section by remarking that $\HF(\phi,+),\HF(\phi,-)$ are
independent of the choices of the local chart $\jmath$ and perturbation
data $H,\delta$.  
They are invariants of the isotopy class of $\phi$ in $\diff_c^+(M)$.


\end{document}